\input amstex\documentstyle {amsppt}  
\pagewidth{12.5 cm}\pageheight{19 cm}\magnification\magstep1
\topmatter
\title Character sheaves on disconnected groups, IV\endtitle
\author G. Lusztig\endauthor
\address Department of Mathematics, M.I.T., Cambridge, MA 02139\endaddress
\thanks Supported in part by the National Science Foundation\endthanks
\endtopmatter   
\document
\define\mpb{\medpagebreak}
\define\co{\eta}

\define\Up{\Upsilon}

\define\bs{\bar s}
\define\bu{\bar u}
\define\bay{\bar y}
\define\baG{\bar\G}
\define\bboc{\bar\boc}

\define\hZ{\hat Z}

\define\uL{\un L}
\define\ufA{\un\fA}

\define\si{\sim}

\define\sqc{\sqcup}

\define\qua{\quad}

\define\hf{\hat f}

\define\tcl{\ti\cl}

\define\baf{\bar f}

\define\bG{\bar G}

\define\bY{\bar Y}

\define\lb{\linebreak}

\define\op{\oplus}

\define\em{\emptyset}

\define\ra{\rangle}

\define\iy{\infty}
\define\m{\mapsto}
\define\do{\dots}
\define\la{\langle}
\define\bsl{\backslash}

\define\lra{\leftrightarrow}

\define\sub{\subset}
\define\bxt{\boxtimes}
\define\T{\times}
\define\ti{\tilde}
\define\nl{\newline}
\redefine\i{^{-1}}

\define\un{\underline}

\define\ot{\otimes}
\define\bbq{\bar{\QQ}_l}

\define\Ad{\text{\rm Ad}}
\define\Hom{\text{\rm Hom}}
\define\End{\text{\rm End}}

\define\Ker{\text{\rm Ker}}

\define\tr{\text{\rm tr}}

\define\a{\alpha}
\redefine\b{\beta}
\redefine\c{\chi}
\define\g{\gamma}
\redefine\d{\delta}
\define\e{\epsilon}
\define\et{\eta}
\define\io{\iota}
\redefine\o{\omega}
\define\p{\pi}
\define\ph{\phi}
\define\ps{\psi}
\define\r{\rho}
\define\s{\sigma}
\redefine\t{\tau}
\define\th{\theta}
\define\k{\kappa}
\redefine\l{\lambda}
\define\z{\zeta}

\redefine\G{\Gamma}

\define\Ps{\Psi}

\define\boc{\bold c}
\define\dd{\bold d}

\define\kk{\bold k}

\define\zz{\bold z}

\define\yy{\bold y}

\redefine\AA{\bold A}

\define\EE{\bold E}
\define\FF{\bold F}

\define\QQ{\bold Q}

\define\VV{\bold V}

\define\ZZ{\bold Z}

\define\ca{\Cal A}

\define\cc{\Cal C}
\define\cd{\Cal D}
\define\ce{\Cal E}
\define\cf{\Cal F}
\define\cg{\Cal G}
\define\ch{\Cal H}
\define\ci{\Cal I}
\define\cj{\Cal J}

\define\cl{\Cal L}

\define\cn{\Cal N}

\define\car{\Cal R}
\define\cs{\Cal S}

\define\cv{\Cal V}
\define\cw{\Cal W}
\define\cz{\Cal Z}
\define\cx{\Cal X}
\define\cy{\Cal Y}

\define\fc{\frak c}

\define\ff{\frak f}

\define\fA{\frak A}

\define\fK{\frak K}

\define\fR{\frak R}

\define\tf{\ti f}
\define\tg{\ti g}

\define\ts{\ti s}

\define\tu{\ti u}

\define\ty{\ti y}

\define\tG{\ti G}

\define\tL{\ti L}

\define\tY{\ti Y}
\define\tZ{\ti Z}

\define\sh{\sharp}

\define\bg{\bar g}

\define\tce{\ti\ce}
\define\bul{\bullet}

\define\uZ{\un Z}
\define\bocd{\boc_\bul}
\define\tiff{\ti\ff}
\define\BOR{B1}
\define\DM{DM}
\define\IC{L2}
\define\ADI{L10}
\define\ADII{L11}
\define\ADIII{L12}
\define\CSV{L13}
\define\SPS{SS}
\define\ST{St}

\head Introduction\endhead
Throughout this paper, $G$ denotes a fixed. not necessarily connected, reductive
algebraic group over an algebraically closed field $\kk$. This paper is a part of a
series (beginning with \cite{\ADI}, \cite{\ADII}, \cite{\ADIII}) which attempts to 
develop a theory of character sheaves on $G$. The numbering of the sections and
references continues that of the earlier Parts.

Assume that $\kk$ is an algebraic closure of a finite field $\FF_q$ and that $G$ has a 
fixed $\FF_q$-rational structure with Frobenius map $F:G@>>>G$. To any triple 
$(L,S,\ce)$ (where $L$ is a Levi of a parabolic of $G^0$, $S$ is an isolated stratum of
the normalizer of $L$, 
with certain properties, and $\ce$ is an irreducible cuspidal local system on $S$) we 
have an associated in 5.6 a (not necessarily irreducible) intersection cohomology 
complex $\fK$ on $G$. If $F(L)=L,F(S)=S$ and we are given an isomorphism 
$F^*\ce@>>>\ce$, there is an induced isomorphism $\ph:F^*\fK@>>>\fK$ hence the 
characteristic function $\c_{\fK,\ph}:G^F@>>>\bbq$ is well defined.

The main result of this paper (Theorem 21.14) is that the functions $\c_{\fK,\ph}$ that
are not identically zero (for various $(L,S,\ce)$ up to $G^{0F}$-conjugacy) form a 
$\bbq$-basis of the vector space $\VV$ of functions $G^F@>>>\bbq$ that are constant on 
$G^{0F}$-conjugacy classes. The proof uses several of the results developed in earlier 
Parts (the generalized Springer correspondence in \S11, the generalized Green functions
in \S15, the character formula in \S16). It also uses the classification of cuspidal 
local systems (this is needed in \S17 which is a preliminary to the proof of Theorem 
21.14).  

A corollary of the main theorem is Theorem 21.21 which states that
the characteristic functions of admissible complexes $A$ such that $F^*A\cong A$ form a
basis for $\VV$. In the connected case such a result was proved in \cite{\CSV} subject
to some mild restrictions on the characteristic. The present proof has no restrictions 
on the characteristic and it makes no use of the orthogonality formulas which will 
appear in a later stage of the theory. 

Another corollary of the main theorem is the construction in \S22 of a "twisted 
induction" map from certain functions on a subgroup of $G^F$ to functions on $G^F$.

Our paper contains also a new characterization of isolated elements (see 18.2) which is
obvious in the connected case but less obvious in the disconnected case.

{\it Notation.} 
If $\G$ is a finite group with a given automorphism $F:\G@>\si>>\G$, the "$F$-twisted 
conjugacy classes" of $\G$ are the orbits of the $\G$-action on $\G$ given by 
$y:w\m F\i(y)wy\i$. 

We shall denote by $\s_G$ the map from $G$ to the set of $G^0$-conjugacy classes of 
quasi-semisimple elements in $G$ defined in 7.1 (where it is denoted by $\s$).

For two elements $a,b$ of a group we set $a^b=b\i ab$. 

Let $p\ge 0$ be the characteristic of $\kk$. 

\mpb

{\it Errata for Part I.}

In 1.11 replace \cite{B, 9.8} by: A.Borel and J.Tits, {\it Groupes r\'eductifs}, 
                                           Publ.Math.IHES 27(1965), 55-150, Lemma 11.1.

6.7: line 3,4: replace $\tce$ by $\p_!\tce$.

p.403 line 5: replace $Z_{L_1}(g_s)$ by $\dim Z_{L_1}(g_s)$.

\head Contents\endhead
17. Properties of cuspidal classes.

18. A property of isolated elements.

19. Properties of cuspidal local systems.

20. Twisted group algebras.

21. Bases.

22. Twisted induction of class functions.

\head 17. Properties of cuspidal classes \endhead
\subhead 17.1\endsubhead
This section contains the proof of a key property (17.13) of "cuspidal conjugacy
classes (see below) which is needed in the proof of the main results of \S21.

\subhead 17.2\endsubhead
A $G^0$-conjugacy class in $G$ is said to be {\it isolated} if one (or equivalently, 
any) element of it is isolated in $G$, see 2.2. We show:

(a) {\it if $G^0$ is semisimple and $g\in G$ is isolated then $g_s$ has finite order.}
\nl
We may assume that $g_u$ is quasi-semisimple in $Z_G(g_s)$. Let $H=Z_G(g_u)$. As in the
proof of Lemma 2.7, we see that $H^0$ is semisimple. By Lemma 2.5, $g_s$ is isolated in
$H$. Replacing $G$ by $H$, we are reduced to the case where $g$ is semisimple. In this 
case we argue by induction on $\dim G$. Assume first that $\dim Z_G(g)<\dim G$. Now 
$Z_G(g)^0$ is semisimple and $g$ is isolated in $Z_G(g)$. By the induction hypothesis
$g$ has finite order. Assume next that $\dim Z_G(g)=\dim G$ that is, $Z_G(g)^0=G^0$. We
can find an integer $n\ge 1$ such that $g^n\in G^0$. Clearly, 
$Z_G(g)^0\sub Z_G(g^n)^0$. Hence $Z_G(g^n)^0=G^0$. Thus, $g^n\in\cz_{G^0}$. Since 
$\cz_{G^0}$ is finite, we see that $g^n$ has finite order. This proves (a).

\subhead 17.3\endsubhead
Let $\boc$ be an isolated $G^0$-conjugacy class in $G$ and let $\cf$ be a local system 
on $\boc$. Let $[\cf]$ be the isomorphism class of $\cf$. We say that $(\boc,\cf)$ (or 
$(\boc,[\cf])$ or $\cf$) is {\it cuspidal} if $\cf$ is $G^0$-equivariant and for any 
proper parabolic $P$ of $G^0$ and any $U_P$-coset $R$ in $N_GP$ we have 
$H_c^d(\boc\cap R,\cf)=0$ where $d$ is $\dim\boc$ minus the dimension of the 
$P/U_P$-conjugacy class of $R/U_P$ in $N_GP/U_P$.

A $G^0$-conjugacy class in $G$ is said to be {\it cuspidal} if it is isolated and if 
it carries some non-zero cuspidal local system.

Let $\boc$ be an isolated $G^0$-conjugacy class in $G$ and let $\cf$ be a 
$G^0$-equivariant local system on $\boc$. Then 

(a) {\it $(\boc,\cf)$ is cuspidal if and only if for some/any $g\in\boc$, and some/any 
unipotent $Z_G(g_s)^0$-conjugacy class $\fc$ of $Z_G(g_s)$ contained in 
$\{u\in Z_G(g_s);u\text{ unipotent, }g_su\in\boc\}$, the local system $j^*\cf$ on $\fc$
is cuspidal relative to $Z_G(g_s)$ (here $j:\fc@>>>\boc$ is $u\m g_su$).}
\nl
The proof of (a) is identical to that of its special case 6.6.

\subhead 17.4\endsubhead
(a) {\it Let $g\in G$ be quasi-semisimple. Then $Z_{G^0}(g)/Z_G(g)^0$ is a
diagonalizable group.}

(b) {\it If in addition, $G^0$ is semisimple, simply connected, we have
$Z_{G^0}(g)=Z_G(g)^0$.} 
\nl
(b) is proved in \cite{\ST, 8.1}; in the closely related case of compact Lie groups, it
goes back to the earlier paper \cite{\BOR, 3.4}. In the case where $G^0$ is semisimple,
(a) is proved in \cite{\ST, 9.1} using (b). For the general case see 
\cite{\DM, 1.6(i)}, \cite{\DM, 1.24}.

\proclaim{Lemma 17.5} Assume that $G$ is such that $G^0$ is semisimple, simply 
connected. Let $s\in G$ be semisimple and $u\in G$ be unipotent such that $su=us$ and 
the $G^0$-conjugacy class of $su$ is cuspidal. Let $z\in\cz_{G^0}$. Let $\uZ=Z_G(s)^0$.
Assume that $g\in G^0,g^s=zg,g^u=yg$ with $y\in\uZ^0$. Then $y=y'{}^uy'{}\i$ for some 
$y'\in\uZ^0$.
\endproclaim
The proof (carried out in 17.6-17.11) consists of a number of steps which reduce us to 
the case where $G^0$ is almost simple, simply connected, in which case we shall use the
classification of cuspidal conjugacy classes (which by 17.3(a) follows from the results
on the classification of unipotent cuspidal conjugacy classes given in \S12 and 
\cite{\IC}). 

\subhead 17.6\endsubhead
In the setup of 17.5 assume that

(a) {\it $G^0$ is almost simple and $G=G^0\T C_m$ (semidirect product, $C_m$ a 
cyclic group of order $m\ge 1$ with generator $e$) with group structure 
$(a,e^t)(a',e^{t'})=(a\a^t(a'),e^{t+t'})$ where $\a:G^0@>>>G^0$ is an 
automorphism of order $m\ge 1$ preserving an \'epinglage and $su$ is of the form 
$(x,e)$ for some $x\in G^0$.}
\nl
If $\cz_{G^0}=\{1\}$ then we must have $z=1$ and we can take $y'=g$. This handles the 
following types of $G^0$:
 
$B_n,C_n$ (with $m=1$, $p=2$), $D_n$ (with $m\in\{1,2\},p=2$)

$D_4$ (with $m=3,p=2$)

$E_6$ (with $m\in\{1,2\},p=3$), $E_7$ (with $p=2$), 

$E_8,F_4$ (with $m\in\{1,2\}$), $G_2$.
\nl
If $\uZ=G^0$ then we must have $z=1$ and we can take $y'=g$. This handles the 
following types of $G^0$:

$A_n$ (with $m=1$), $A_n$ (with $m=2,p=2$),

$E_6$ (with $m=1,p=2$, $\uZ=G^0$), $E_7$ (with $p=3$, $\uZ=G^0$).
\nl
Let $\boc_1$ (resp. $\boc_2$) be the $\uZ^0$-conjugacy class of $u$ (resp. $uy$) in 
$Z_G(s)$. Then $\Ad(g):Z_G(s)@>\si>>Z_G(s)$ carries $\boc_1$ to $\boc_2$. Thus 
$\boc_1,\boc_2$ are two cuspidal unipotent $\uZ$-conjugacy classes in the same 
connected component of $Z_G(s)$. It is enough to show that $\boc_1=\boc_2$. (Then 
$uy=y'uy'{}\i$ for some $y'\in\uZ$.) If there is only one cuspidal unipotent
$\uZ$-conjugacy class in $u\uZ$ then clearly $\boc_1=\boc_2$. This handles the 
following types of $G^0$:

$A_n$ (with $m=2,p\ne 2$), $C_n$, $D_4$ (with $m=3,p\ne 2$),

$E_6$ (with $m=1,p\ne 2,3$), $E_6$ (with $m=1,p=2,\uZ\ne G^0$), 

$E_6$ (with $m=2,p\ne 3$),

$E_7$ (with $p\ne 2,3$), $E_7$ (with $p=3,\uZ\ne G^0$).
\nl
The cases not covered by the arguments above are with $G_0$ of type $B_n$ or $D_n$ with
$m\in\{1,2\}$ and $p\ne 2$. In each of these cases there are at most two cuspidal
unipotent $\uZ$-conjugacy classes in $Z_G(s)$ (they are automatically contained in
$\uZ$). Each of these classes is stable under any automorphism of $\uZ$; in particular,
under $\Ad(g)$. Hence these cases are settled. Thus, Lemma 17.5 holds in the present 
case.

\subhead 17.7\endsubhead
In the setup of 17.5 assume that 

(a) {\it $G^0$ is almost simple and $G=G^0\T C_n$ (semidirect 
product, $C_n$ a cyclic group of order $n\ge 1$ with generator $e$) with group 
structure $(a,e^t)(a',e^{t'})=(a\a^t(a'),e^{t+t'})$ where $\a:G^0@>>>G^0$ is an 
automorphism preserving an \'epinglage such that $\a^n=1$ and $su$ is of the form 
$(x,e)$ for some $x\in G^0$.}
\nl
Let $m$ be the order of $\a$. Thus,
$n/m\in\ZZ$. Let $\bG=G^0\T C_m$ (semidirect product, $C_n$ a cyclic group of order
$m\ge 1$ with generator $e'$) with group structure 
$(a,e'{}^t)(a',e'{}^{t'})=(a\a^t(a'),e'{}^{t+t'})$. Let $\p:G@>>>\bG$,
$(g,e^t)\m(g,e'{}^t)$. Then $\p$ induces $G^0@>\si>>\bG^0$ with kernel
$K=\{1,e^m,e^{2m},\do\}$. From the definitions we see that the $\bG^0$-conjugacy class 
of $\p(s)\p(u)$ is cuspidal. Applying 17.6 to $\bG,\p(s),\p(u),\p(g),\p(z),\p(y)$
instead of $G,s,u,g,z,y$, we see that there exists $y'\in G^0$ such that $sy'=y'sk$,
$uy=y'uy'{}\i k'$ with $k,k'\in K$. Applying the homomorphism 
$\r:G@>>>C_n,(g,e^t)\m e^t$ we get 

$\r(s)\r(y')=\r(y')\r(s)\r(k)$, $\r(u)\r(y)=\r(y')\r(u)\r(y')\i\r(k')$, 

$\r(u)\i\r(g)\r(u)=\r(y)\r(g)$.
\nl
Using the commutativity of $C_n$ we deduce $\r(y)=1,\r(k')=1,\r(k)=1$. Since 
$\r:K@>>>C_n$ is injective it follows that $k=k'=1$. Thus, Lemma 17.5 holds in the 
present case.

\subhead 17.8\endsubhead
In the setup of 17.5 assume that 

(a) {\it $G^0$ is almost simple and $G$ is generated by the connected component $D$
that contains $su$.}
\nl
We can find $d\in D$ such that $\Ad(d):G^0@>>>G^0$ preserves an \'epinglage. Then 
$d$ has order $n<\iy$. Let $\tG=G^0\T C_n$ (semidirect product, $C_n$ a cyclic group 
of order $n\ge 1$ with generator $e$) with group structure 
$(a,e^t)(a',e^{t'})=(a\Ad(d)^t(a'),e^{t+t'})$. Now $\p:\tG@>>>G,\p(g,e^t)=gd^t$ is a 
group homomorphism with kernel $K=\{(d^t,e^{-t});d^t\in G^0\}$; it induces an
isomorphism $\tG^0@>\si>>G^0$. Let $x\in\tG$ be s.t. $\p(x)=su$, $x$ of the form 
$(x_0,e),x_0\in G^0$. We have $(g,1)x_s(g\i,1)=x_s(z,1)k$, with $k\in K$. Taking images
in $\tG/G^0$ we see that $k$ goes to the neutral element hence $k\in G^0$. But 
$K\cap G^0=\{1\}$ so that $k=1$. We have $(g,1)x_u(g\i,1)=x_u(y,1)k'$ with $k'\in K$. 
As above we see that $k'=1$. We have $(y,1)x_s(y\i,1)=x_sk''$ with $k''\in K$. As above
we see that $k''=1$. From the definitions we see that the $\tG^0$-conjugacy class of 
$x$ is cuspidal. Applying 17.7 to $\tG,x_s,x_u,(g,1),(y,1),(z,1)$ instead of 
$G,s,u,g,y,z$ we find $y'\in G^0$ such that $(y',1)x_s(y'{}\i,1)=x_s$, 
$(y,1)=x_u\i(y',1)x_u(y'{}\i,1)$. Applying $\p$ we get $y'sy'{}\i=s,y=u\i y'uy'{}\i$. 
Thus Lemma 17.5 holds in the present case.

\subhead 17.9\endsubhead
In the setup of 17.5 assume that $G$ has no closed connected normal subgroup other than
$G^0$ and $\{1\}$. Assume also that $p=0$. We have $G^0=\prod_{j\in\ZZ/b}H_j$ where 
$H_j$ is connected, simply connected, almost simple, $b\ge 1$ and $sH_js\i=H_{j+1}$, 
$u\in G^0$. Set $g=(g_j),z=(z_j),y=(y_j),u=(u_j)$ where $g_j\in H_j,z_j\in\cz_{H_j}$,
$y_j\in H_j,u_j\in H_j$. We have $g_{j+1}^s=z_jg_j$, $g_j^{u_j}=y_jg_j$, 
$y_{j+1}^s=y_j$, $z_j^u=z_j$, $u_{j+1}^s=u_j$. Let $G'$ be the subgroup of $G$ 
generated by $H_0,s^b$. This is a closed subgroup with identity component $H_0$ since 
$s$ has finite order (see 17.2; recall that $su$ is isolated in $G$). We have 
$g_0^{u_0}=y_0g_0,g_0^{s^b}=\zz g_0$, where 
$\zz=z_{b-1}^{s^{b-1}}\do z_1^sz_0\in\cz_{H_0}$. We have $y_0^{s^b}=y_0$. Also 
$s^bu_0=u_0s^b$. We show that the $H_0$-conjugacy class of $s^bu_0$ is cuspidal. By 
17.3(a), it is enough to show that:

(i) the $Z_{H_0}(s^b)^0$-conjugacy class of $u_0$ is cuspidal in $Z_{G'}(s^b)^0$,

(ii) $s^bu_0$ is isolated in $G'$.
\nl
Since the $G^0$-conjugacy class of $su$ is cuspidal, we see from 17.3(a) that:

(iii) the $Z_G(s)^0$-conjugacy class of $u$ is cuspidal in $Z_G(s)$,

(iv) $su$ is isolated in $G$.
\nl
Now $Z_G(s)^0=Z_{G^0}(s)$ (see 17.4(b)) consists of all $(x_j)$ where $x_j\in H_j$ 
satisfy $x_{j+1}^s=x_j$. We may identify $Z_{G^0}(s)=Z_{H_0}(s^b)$ and (i) follows. We 
prove (ii). From (iv) we see that $\cz_{Z_{G^0}(s)}\cap Z_{G^0}(u)$ is finite. Hence 
$f\in Z_{H_0}(s^b)$ subject to $f^{u_0}=f$ has finitely many possible values. Hence 
(iv) holds.

Applying 17.8 to $G',s^b,u_0,g_0,y_0,\zz$ instead of $G,s,u,g,y,z$ we find $\ty\in H_0$
such that $\ty^{s^b}=\ty$, $y_0=\ty^{u_0}\ty\i$. Set $y'_j=\ty^{s^{-j}}\in H_j$. 
Clearly, $y'_{j+1}{}^s=y'_j,y_0y'_0=y'_0{}^{u_0}$. Hence $y_jy'_j=y'_j{}^{u_j}$. (We
have 

$y_jy'_j=(y_0y'_0)^{s^{-j}}=(y'_0{}^{u_0})^{s^{-j}}
=(y'_0{}^{s^{-j}})^{u_j}=y'_j{}^{u_j}$.) 
\nl
Hence setting $y'=(y'_j)$ we have $y'{}^s=y'$, $y=y'{}^uy'{}\i$. Thus Lemma 17.5 holds 
in the present case.

\subhead 17.10\endsubhead
In the setup of 17.5 assume that 

(a) {\it $G$ has no closed connected normal subgroup other than $G^0,\{1\}$ and $p>1$.}
\nl
We have $G^0=\prod_{i\in\ZZ/a,j\in\ZZ/b}H_{ij}$ where $H_{ij}$ is connected, simply 
connected, almost simple, $a\ge 1,b\ge 1$ and $uH_{ij}u\i=H_{i+1,j}$, 
$sH_{ij}s\i=H_{i,j+1}$. Let $G'$ be the subgroup of $G$ generated by $H_{00}$, 
$u^a,s^b$. This is a closed subgroup with identity component $H_{00}$ since $s$ has 
finite order (see 17.2; recall that $su$ is isolated) and $u$ has finite order, power 
of $p$. Now $a$ is a power of $p$ and $b$ is prime to $p$. Set 
$g=(g_{ij}),z=(z_{ij}),y=(y_{ij})$ where $g_{ij}\in H_{ij}$, $z_{ij}\in\cz_{H_{ij}}$, 
$y_{ij}\in H_{ij}$. We have 

$g_{i,j+1}^s=z_{ij}g_{ij}$, $g_{i+1,j}^u=y_{ij}g_{ij}$, $y_{i,j+1}^s=y_{ij}$, 
$z_{i+1,j}^u=z_{ij}$.
\nl
(The last equation follows from $uz=zu$: we have 

$z^ug^u=(g^s)^u=(g^u)^s=(yg)^s=yzg=zyg=zg^u$ 
\nl
hence $z^u=z$.) We have

$g_{00}^{u^a}=\yy g_{00},\yy=y_{a-1,0}^{u^{a-1}}\do y_{1,0}^uy_{0,0}\in H_{00}$,

$g_{00}^{s^b}=\zz g_{00},\zz=z_{0,b-1}^{s^{b-1}}\do z_{0,1}^sz_{0,0}\in\cz_{H_{00}}$.
\nl
Since $y_{ij}^{s^b}=y_{ij}$ we have $\yy^{s^b}=\yy$. Also, $s^bu^a=u^as^b$. We show 
that the $H_{00}$-conjugacy class of $s^bu^a$ is cuspidal. By 17.3(a) it is enough to 
show that:

(i) the $Z_{H_{00}}(s^b)^0$-conjugacy class of $u^a$ is cuspidal in $Z_{G'}(s^b)^0$,

(ii) $s^bu^a$ is isolated in $G'$.
\nl
Since the $G^0$-conjugacy class of $su$ is cuspidal, we see from 17.3(a) that:

(iii) the $Z_G(s)^0$-conjugacy class of $u$ is cuspidal in $Z_G(s)$,

(iv) $su$ is isolated in $G$.
\nl
Now $Z_G(s)^0=Z_{G^0}(s)$ (see 17.4(b)) consists of all $(x_{ij})$ where 
$x_{ij}\in H_{ij}$ satisfy $x_{i,j+1}^s=x_{ij}$. For each $i\in\ZZ/a$ let 
$F_i=\prod_{j\in\ZZ/b}H_{ij}$. Then 

$sF_is\i=F_i$, $uF_iu\i=F_{i+1}$, $Z_{G^0}(s)=\prod_{i\in\ZZ/a}Z_{F_i}(s)$.
\nl
By an argument in 12.5(e) applied to 
$Z_G(s),Z_G(s)^0=\prod_iZ_{F_i}(s),u$ instead of $G,G^0=\prod_iH_i,u$, we see that the 
$Z_{F_0}(s)$-conjugacy class of $u^a$ is cuspidal in the subgroup generated by 
$Z_{F_0}(s),u^a$. We may identify $Z_{F_0}(s)=Z_{H_{00}}(s^b)$ and (i) follows. We 
prove (ii). From (iv) we see that $\cz_{Z_{G^0}(s)}\cap Z_{G^0}(u)$ is finite. Hence if
$(f_i)$ satisfies $f_i\in\cz_{Z_{F_i}(s)},f_{i+1}^u=f_i$ then $f_i$ has finitely many 
values. Hence $f_0\in\cz_{Z_{F_0}(s)}$, subject to $f_0^{u^a}=f_0$ has finitely many 
values. Hence $f_0\in Z_{H_{00}}(s^b)$ subject to $f_0^{u^a}=f_0$ has finitely many 
values. Hence (ii) holds.

Applying 17.8 to $G',s^b,u^a,g_{00},\yy,\zz$ instead of $G,s,u,g,y,z$ we find 
$\ty\in H_{00}$ such that $\ty^{s^b}=\ty$, $\yy=\ty^{u^a}\ty\i$. Set
$y'_{ij}\in H_{ij}$ by $y'_{0j}=\ty^{s^{-j}}$ for $j\in\ZZ/b$,

$y'_{ij}=y_{i-1,0}^{u\i s^{-j}}\do y_{1,0}^{u^{-i+1}s^{-j}}y_{0,0}^{u^{-i}s^{-j}}
\ty^{u^{-i}s^{-j}}$
\nl
for $i=1,\do,a-1$ and $j\in\ZZ/b$. Clearly, $y'_{i,j+1}{}^s=y'_{ij}$. Moreover,
$y_{ij}y'_{ij}=y'_{i+1,j}{}^u$ for $i=0,1,\do,a-2$. The same holds for $i=a-1$:
$$\align&y_{a-1,j}y'_{a-1,j}=y_{a-1,0}^{s^{-j}}y_{a-2,0}^{u\i s^{-j}}\do 
y_{1,0}^{u^{-a+2}s^{-j}}y_{0,0}^{u^{-a+1}s^{-j}}\ty^{u^{-a+1}s^{-j}}  \\&
=\yy^{u^{-a+1}s^{-j}}\ty^{u^{-a+1}s^{-j}} =\ty^{us^{-j}}=y'_{0j}{}^u=y'_{aj}{}^u.
\endalign$$
Hence setting $y'=(y'_{ij})\in G^0$ we have $y'{}^s=y',y=y'{}^uy'{}\i$. Thus Lemma 17.5
holds in the present case.

\subhead 17.11\endsubhead
We now prove Lemma 17.5 by induction on $\dim G$. If $\dim G=0$ the result is trivial. 
We now assume that $\dim G>0$. Assume first that $G^0=G_1\T G_2$ where $G_i\ne\{1\}$ 
are connected, simply connected, normal in $G$. Let 
$G'_1=G/G_2,G'_2=G/G_1,G'=G'_1\T G'_2$. Then $G\sub G',G^0=G'{}^0$. We have 

$s=(s_1,s_2),u=(u_1,u_2),z=(z_1,z_2),g=(g_1,g_2),y=(y_1,y_2)$
\nl
where $s_i$ is semisimple
in $G'_i$, $u_i$ is unipotent in $G'_i$, $z_i\in\cz_{G_i}$, $g_i\in G_i$,
$y_i\in Z_{G'_i}(s_i)^0$. We have $s_iu_i=u_is_i$, $g_i^{s_i}=z_ig_i$,
$g_i^{u_i}=y_ig_i$. Also the $G_i$-conjugacy class of $s_iu_i$ is cuspidal. By the
induction hypothesis, we can find $y'_i\in Z_{G'_i}(s_i)^0$ with 
$y_i=y'_i{}^{u_i}y'_i{}\i$. Let $y'=(y'_1,y'_2)$. Then $y'\in Z_{G^0}(s)$, 
$y=y'{}^uy'{}\i$. Thus Lemma 17.5 holds in the present case.

Next we assume that no decomposition $G^0=G_1\T G_2$ as above exists. Then the result
follows from 17.9, 17.10. The lemma is proved.

\proclaim{Lemma 17.12} Let $s\in G$ be semisimple and $u\in G$ be unipotent such 
that $su=us$ and the $G^0$-conjugacy class of $su$ is cuspidal. Assume that $g\in G^0$,
$g^s=g,g^u=yg$ with $y\in Z_G(s)^0$. Then $y=y'{}^uy'{}\i$ for some $y'\in Z_G(s)^0$.
\endproclaim
Assume first that $G^0$ is semisimple and that 

(a) {\it there exists an element $d$ in the connected component of $G$ that contains 
$su$ such that $\Ad(d):G^0@>>>G^0$ preserves an \'epinglage and such that 
$\{d^t;t\in\ZZ\}\cap G^0=\{1\}$.}
\nl
Using \cite{\ST, 9.16} we can find a reductive group $\tG$ with $\tG$ semisimple, 
simply connected and a surjective homomorphism of algebraic groups $\p:\tG@>>>G$ such 
that $K=\Ker\p\sub\cz_{\tG^0}$. Pick $\ts'\in\tG$ semisimple, $\tu\in\tG$ unipotent 
such that $\p(\ts')=s,\p(\tu)=u$. Then $\tu\ts'=\ts'\tu k$ with $k\in K$. Now 
$\tu k=\ts'{}\i\tu\ts'$ is unipotent. Since $\tu$ normalizes $K$ (a diagonalizable 
group) it follows that $\tu k=k'\tu k'{}\i$ for some 
$k'\in K$. Set $\ts=\ts'k'$. Then $\ts$ is semisimple (since $\ts'$ normalizes $K$, a 
diagonalizable group) and $\p(\ts)=s$, $\tu\ts=\ts\tu$. Since $\p(\ts\tu)=su$ we see 
that the $G^0$-conjugacy class of $su$ is cuspidal. Let $g'\in\tG^0$ be such that 
$\p(g')=g$. Since $\hZ=\{x\in\tG;\ts x\in x\ts\cz_{\tG^0}\}$ has identity component 
$Z_{\tG}(\ts)^0$, we see that $Z_{\tG}(\ts)^0@>>>Z_G(s)^0$ is surjective. Hence we can 
find $\ty\in Z_{\tG}(\ts)^0$ such that $\p(\ty)=y$. We have $\p(g'{}^{\tu})=\p(\ty g')$
hence $g'{}^{\tu}=k'\ty g'$ for some $k'\in\G$. Hence $g'\tu g'{}\i=\tu k'\ty$, 
equality in $\hZ$. Hence $\tu k'$ is unipotent in $\hZ/\hZ^0$. Since the image of 
$\tu k'$ in $\hZ/\hZ^0$ normalizes the image of $K$ in $\hZ/\hZ^0$, we see that there 
exists $k_1\in K$ such that $\tu k'=k_1\tu k_1\i$ (equality in $\hZ/\hZ^0$) hence
$\tu k'=k_1\tu k_1\i\ty_1$ (equality in $\hZ$) for some $\ty_1\in Z_{\tG}(\ts)^0$.
Hence $g'\tu g'{}\i=k_1\tu k_1\i\ty_1\ty=k_1\tu\ty''k_1\i$ for some 
$\ty''\in Z_{\tG}(\ts)^0$. Set $\tg=k_1\i g'$. Then $\tg\in\tG^0$, 
$\tg^{\tu}=\ty''\tg$, $\p(\tg)=g$. We have $\p(\tg^{\ts})=\p(\tg)$ hence 
$\tg^{\ts}=z\tg$ for some $z\in K$. Applying Lemma 17.5 to $\tG,\ts,\tu,\tg,\ty'',z$ 
instead of $G,s,u,g,y,z$ we see that $\ty''=\ty'{}^{\tu}\ty'{}\i$ for some 
$\ty'\in Z_{\tG}(\ts)^0$. Let $y'=\p(\ty')$. Then $\p(\ty'')=y'{}^uy'{}\i$,
$y'\in Z_G(s)^0$. Also, $g^u=\p(\ty'')g,g^u=yg$ hence $\p(\ty'')=y$ and
$y=y'{}^uy'{}\i$. Thus the lemma holds in the present case.

Next assume that $G^0$ is semisimple and $G$ is generated by a connected component 
$D$. We can find $d\in D$ such that $\Ad(d):G^0@>>>G^0$ preserves an \'epinglage. Then 
$d$ has order $n<\iy$. Let $G'=G^0\T C_n$ (semidirect product, $C_n$ a cyclic group of
order $n\ge 1$ with generator $e$) with group structure 
$(a,e^t)(a',e^{t'})=(a\Ad(d)^t(a'),e^{t+t'})$. Now $\p':G'@>>>G,\p(g,e^t)=gd^t$ is a 
group homomorphism with kernel $K'=\{(d^t,e^{-t});d^t\in\cz_{G^0}\}$; it induces an
isomorphism $G'{}^0@>\si>>G^0$. Then $G'$ is as in the first part of the proof.
Let $x\in G'$ be such that $\p'(x)=su$. From the definitions we see that the
$G'{}^0$-conjugacy class of $x$ is cuspidal. For any $h\in Z_G(s)$ we have 
$(h,1)x_s(h\i,1)=x_sk$ with $k\in K'$. Taking images in $G'/G'{}^0$ we see that $k$ 
goes to $1$ hence $k\in G'{}^0\cap K'=\{1\}$ hence $k=1$. Thus, 
$(h,1)\in Z_{G'}(x_s)$. It follows that for any $h\in Z_G(s)^0$ we have 
$(h,1)\in Z_{G'}(x_s)^0$. In particular $(y,1)\in Z_{G'}(x_s)^0$. In the same way we
see that $(g,1)x_s(g\i,1)=x_s$, $(g,1)x_u(g\i,1)=x_u(y,1)$ (compare 17.8). From the
first part of the proof we see that there exists $(y',1)\in Z_{G'}(x_s)^0$ such that 
$(y,1)=x_u\i(y',1)x_u(y'{}\i,1)$. It follows that $y'\in Z_G(s)^0$, $y=u\i y'uy'{}\i$.
Thus the lemma holds in the present case.

Next, assume that $G^0$ is semisimple (but there is no assumption on $G/G^0$). Let 
$G_1$ be the subgroup of $G$ generated by the connected component that contains $su$. 
By the earlier part of the proof, the lemma holds for $G_1,s,u,g,y$ instead of 
$G,s,u,g,y$. But then it automatically holds for $G,s,u,g,y$.

Finally, we consider the general case. Let $\p'':G@>>>\bG=G/\cz_{G^0}^0$ be the obvious
homomorphism. Let $\bs=\p''(s),\bu=\p''(u)$. Then the $\bG^0$-conjugacy class of 
$\bs\bu=\bu\bs$ is cuspidal. Let $\bg=\p''(g)\in\bG^0,\bay=\p''(y)$. Then
$\bg^{\bs}=\bg,\bg^{\bu}=\bay\bg$, $\bay\in Z_{\bG}(\bs)^0$. Since the lemma holds for 
$\bG$ instead of $G$, we have $\bay=\bay'{}^{\bu}\bay'{}\i$ for some 
$\bay'\in Z_{\bG}(\bs)^0$. Let $\tZ=\{x\in G;xsx\i\in s\cz_{G^0}^0\}$. Then $\p''$ 
induces a surjective homomorphism $\tZ@>>>Z_{\bG}(\bs)$. Moreover, $\tZ^0=Z_G(s)^0$ 
hence $\p''$ 
induces a surjective homomorphism $Z_G(s)^0@>>>Z_{\bG}(\bs)^0$. Hence we can find 
$y'_1\in Z_G(s)^0$ such that $\p''(y'_1)=\bay'$. We have $y=y'_1{}^uy'_1{}\i z$ for 
some $z\in\cz_{G^0}^0$. Thus $uy=y'_1uy'_1{}\i z$. Since $uy=gug\i$ is unipotent, we 
see that $y'_1uy'_1{}\i z$ is unipotent. Also, $z=y'_1(u\i y'_1{}\i u)y\in Z_G(s)^0$ 
(since $y,y'_1\in Z_G(s)^0,u\in Z_G(s)$) hence $z\in Z_G(s)^0\cap\cz_{G^0}^0$. 

Assume first that $p=0$. We set $y'=y'_1$. Then $y'uy'{}\i$ being unipotent is in $G^0$
hence it commutes with $z$. Since $y'uy'{}\i z$ is unipotent, we have $z=1$ and 
$y=y'{}^uy'{}\i$, as required.

Next assume that $p>1$. Then $y'_1uy'_1{}\i$ has finite order and, being in $Z_G(s)$, 
it normalizes $H=Z_G(s)^0\cap\cz_{G^0}^0$. Hence, if $H'$ is the subgroup of $G$ 
generated by $H$ and $y'_1uy'_1{}\i$, we see that $H'$ contains $H$ as a normal 
subgroup of finite index, a power of $p$. Since $H$ is diagonalizable, it follows that 
any two unipotent elements of $H'$ in the same $H$-coset are $H$-conjugate. In 
particular, the unipotent elements $uy=y'_1uy'_1{}\i z,y'_1uy'_1{}\i$ of $H'$ are 
$H$-conjugate. Hence $uy=\z y'_1uy'_1{}\i\z\i$ for some $\z\in H$. We set $y'=\z y'_1$.
Then $y'\in Z_G(s)^0$ and $uy=y'uy'{}\i$. The lemma is proved.

\proclaim{Proposition 17.13} Let $\boc$ be a cuspidal $G^0$-conjugacy class in $G$. Let
$s\in G$ be the semisimple part of some element of $\boc$ and let $\d$ be a connected 
component of $Z_G(s)$. Then $\{u\in Z_G(s)\text{ unipotent },u\in\d,su\in\boc\}$ is a 
single $Z_G(s)^0$-conjugacy class.
\endproclaim
We must show that, if $u,u'\in Z_G(s)$ are unipotent, $u'\in uZ_G(s)^0$ and 
$su\in\boc,su'\in\boc$, then $u,u'$ are $Z_G(s)^0$-conjugate. We can find $g\in G^0$ 
such that $gsug\i=su'$. Then $gsg\i=s,gug\i=u'=uy$ where $y\in Z_G(s)^0$. Hence 
$g^s=g,g^u=yg$. By Lemma 17.12 we can find $y'\in Z_G(s)^0$ such that 
$y=u\i y'uy'{}\i$. Then $u'=y'uy'{}\i$. The proposition is proved. 

\head 18. A property of isolated elements\endhead
\subhead 18.1\endsubhead
This section contains a characterization of isolated elements of $G$.
(This will not be used later in this Part.)

\proclaim{Proposition 18.2} Let $s\in G$ be semisimple and $u\in G$ be unipotent such 
that $su=us$. Then $su$ is isolated in $G$ if and only if $s$ is isolated in $G$.
Equivalently, we have $(\cz_{Z_G(s)^0}\cap Z_G(u))^0=(\cz_{G^0}\cap Z_G(su))^0$ if and 
only if $\cz_{Z_G(s)^0}^0=(\cz_{G^0}\cap Z_G(s))^0$.
\endproclaim
Clearly, if $s$ is isolated in $G$ then $su$ is isolated in $G$. The proof of the 
converse is given in 18.3-18.12. 

{\it In 18.3-18.12, it is assumed that $su$ is isolated in $G$.}

\subhead 18.3\endsubhead
In the setup of 18.2 assume that $u\in G^0$. (This condition is automatically satisfied
if $p=0$.) By assumption, we have $u\in Z_{G^0}(s)$. The image of $u$ in 
$Z_{G^0}(s)/Z_G(s)^0$ is semisimple (by 17.4(a)) and unipotent hence is $1$. Thus, 
$u\in Z_G(s)^0$. It follows that $(\cz_{Z_G(s)^0}\cap Z_G(u))^0=\cz_{Z_G(s)^0}^0$ and
$(\cz_{G^0}\cap Z_G(su))^0=(\cz_{G^0}\cap Z_G(s))^0$. Hence the condition 
$(\cz_{Z_G(s)^0}\cap Z_G(u))^0=(\cz_{G^0}\cap Z_G(su))^0$ implies that 
$\cz_{Z_G(s)^0}^0=(\cz_{G^0}\cap Z_G(s))^0$. Thus $s$ is isolated in $G$.

\subhead 18.4\endsubhead
In the setup of 18.2 assume that $G^0$ is semisimple, that $s\in G^0$. We show that 
$s$ is isolated in $G$. We may assume that $u$ is unipotent, quasi-semisimple in 
$Z_G(s)$. Let $\cz=\cz_{Z_G(s)^0}^0$. Assume that $\cz\ne\{1\}$. Let $L=Z_{G^0}\cz$, a 
Levi of a parabolic of $G^0$. Since $s\in G^0$, we have $\cz_L^0=\cz$. Since $su$ is 
quasi-semisimple, we can find a Borel $B$ and a maximal torus $T$ of $B$ that are 
normalized by $su$. Since $s\in G^0$, we have $s\in T$. Hence $T\sub Z_G(s)^0$ so that 
$\cz\sub T$ and $T\sub L$. Let $\b=B\cap L$, a Borel of $L$. Let 
$\Pi\sub V=\Hom(T,\kk^*)\ot\QQ$ be the set of simple roots of $G^0$ with respect to 
$T,B$ (in particular, the corresponding root subgroups are contained in $U_B$). Let $Q$
be the basis of $V^*=\Hom(\kk^*,T)\ot\QQ$ dual to $\Pi$. There is a unique subset $Q_1$
of $Q$ which is a basis for the subspace $V^*_1=\Hom(\kk^*,\cz)\ot\QQ$ of $V^*$. Since
$\cz\ne\{1\}$ we have $V^*_1\ne 0$ hence $Q_1\ne\em$. Now $u$ normalizes $Z_G(s)^0$ 
hence $u\cz u\i=\cz,uLu\i=u,u\b u\i=\b$. Hence the automorphism of $V^*_1$ induced by
$\Ad(u)$ preserves $Q$ and $Q_1$ and the sum of elements in $Q_1$ is a non-zero 
$\Ad(u)$-invariant vector. Thus, $\Ad(u):V_1^*@>>>V_1^*$ has a non-zero fixed point 
set. It follows that $\dim(\cz\cap Z_G(u))>0$ contradicting the assumption that $su$ is
isolated in $G$. We have proved that $\cz=\{1\}$. Hence $s$ is isolated in $G$.

\subhead 18.5\endsubhead
In the setup of 18.2 assume that $G^0$ is semisimple, simply connected and 17.6(a) 
holds. Since $m$ in 17.6(a) is $1$ or a prime number, we have three cases:

(i) $u\in G^0,s\notin G^0$ and $1<m\ne p$;

(ii) $u\notin G^0,s\in G^0$ and $1<m=p$;

(iii) $u\in G^0,s\in G^0$ and $1=m$.
\nl
In cases (ii),(iii) we have $s\in G^0$ hence by the argument in 18.4 we see that $s$ is
isolated in $G$. In cases (i),(iii) we have $u\in G^0$ hence by the argument in 18.3 we
see that $s$ is isolated in $G$. 

\subhead 18.6\endsubhead
In the setup of 18.2 assume that $G^0$ is semisimple, simply connected and 17.7(a) 
holds. Let $\p:G@>>>\bG$ be as in 17.7. Let $\bs=\p(s),\bu=\p(u)$. One checks that 
$\bs\bu$ is isolated in $\bG$. Applying 18.5 to $\bG,\bs,\bu$ instead of $G,s,u$ we see
that $\cz_{Z_{\bG}(\bs)^0}^0=\{1\}$. Hence $\cz_{Z_G(s)^0}^0=\{1\}$. Thus $s$ is 
isolated in $G$.

\subhead 18.7\endsubhead
In the setup of 18.2 assume that $G^0$ is semisimple, simply connected and 17.8(a) 
holds. Let $\p:\tG@>>>G,x\in\tG$ be as in 17.8. One checks that $x$ is isolated in
$\tG$. Applying 18.6 to $\tG,x_s,x_u$ instead of $G,x,s$ we see that 
$\cz_{Z_{\tG}(\ts)^0}^0=\{1\}$. Hence $\cz_{Z_G(s)^0}^0=\{1\}$. Thus $s$ is isolated 
in $G$.

\subhead 18.8\endsubhead
In the setup of 18.2 assume that $G^0$ is semisimple, simply connected and 17.10(a) 
holds. Let $a,b,H_{ij},G',F_k$ be as in 17.10. As in 17.10 we see that $s^bu^a$ is 
isolated in $G'$. Applying 18.7 to $G',s^b,u^a$ instead of $G,s,u$ we see that 
$\cz_{Z_{G'}(s^b)^0}^0=\{1\}$ that is $\cz_{Z_{H_{00}}(s^b)^0}^0=\{1\}$. Hence 
$\cz_{Z_{F_0}(s)^0}^0=\{1\}$. Now $u^kF_0u^{-k}=F_k$ and since $us=su$, $\Ad(u^k)$ is 
an isomorphism $Z_{F_0}(s)^0@>\si>>Z_{F_k}(s)^0$ hence we have
$\cz_{Z_{F_k}(s)^0}^0\cong\cz_{Z_{F_0}(s)^0}^0=\{1\}$. Hence
$\cz_{Z_G(s)^0}^0=\prod_i\cz_{Z_{F_i}(s)^0}^0=\{1\}$. Thus $s$ is isolated in $G$.

\subhead 18.9\endsubhead
In the setup of 18.2 assume that $G^0$ is semisimple, simply connected. We show that 
$s$ is isolated in $G$ by induction on $\dim G$. If $\dim G=0$ the result is trivial. 
We now assume that $\dim G>0$. Assume first that $G^0=G_1\T G_2$ where $G_i\ne\{1\}$ 
are connected, simply connected, normal in $G$. Let $G'_i,s_i,u_i$ be as in 17.11.
Then $s_iu_i=u_is_i$ is isolated in $G'_i$. By the induction hypothesis we have
$\cz_{Z_{G_i}(s_i)^0}^0=\{1\}$ for $i=1,2$. Now
$Z_G(s)^0=Z_{G_1}(s_1)^0\T Z_{G_2}(s_2)^0$ hence 
$\cz_{Z_G(s)^0}^0=\cz_{Z_{G_1}(s_1)^0}^0\T\cz_{Z_{G_2}(s_2)^0}^0=\{1\}$.
Thus $s$ is isolated in $G$.

Next we assume that no decomposition $G^0=G_1\T G_2$ as above exists. If $p>1$, then 
18.8 shows that $s$ is isolated in $G$. If $p=0$ then 18.3 shows that $s$ is isolated 
in $G$. This completes the inductive proof.

\subhead 18.10\endsubhead
In the setup of 18.2 assume that $G^0$ is semisimple and that 17.12(a) holds.
Let $\p:\tG@>>>G,\ts,\tu$ be as in 17.12. We show that $\ts\tu$ is isolated in $\tG$.
Let $x\in(\cz_{Z_{\tG}(\ts)^0}^0\cap Z_{\tG}(\tu))^0$. Then 
$\p(x)\in\cz_{Z_G(s)^0}^0\cap Z_G(u)$. Indeed, the map $\p:Z_{\tG}(\ts)^0@>>>Z_G(s)^0$ 
is a surjective, finite covering of connected reductive groups hence it restricts to a 
surjective map

($*$) $\cz_{Z_{\tG}(\ts)^0}^0@>>>\cz_{Z_G(s)^0}^0$.
\nl
We see that $\p(x)\in(\cz_{Z_G(s)^0}^0\cap Z_G(u))^0$. Hence 
$\p(x)\in(\cz_{G^0}^0\cap Z_G(u))^0$ and $\p(x)=1$. Hence $x\in\Ker\p$, a finite group.
It follows that $x=1$ and $\ts\tu$ is isolated in $\tG$. Applying 18.9 to $\tG,\ts,\tu$
instead of $G,s,u$, we see that $\cz_{Z_{\tG}(\ts)^0}^0=\{1\}$. Let 
$y\in(\cz_{Z_G(s)^0})^0$. By the surjectivity of ($*$) we have $y=\p(y')$ where 
$y'\in\cz_{Z_{\tG}(\ts)^0}^0$. Hence $y'=1$ and $y=1$. We see that 
$(\cz_{Z_G(s)^0})^0=\{1\}$. Thus $s$ is isolated in $G$.

Next we assume, in the setup of 18.2 that $G^0$ is semisimple and $G/G^0$ is cyclic.
Let $\p':G'@>>>G,x$ be as in 17.12. Then $x$ is isolated in $G'$. By the argument above
applied to $G'$ instead of $G$ we see that $x_s$ is isolated in $G'$. It follows 
immediately that $s$ is isolated in $G$.

\subhead 18.11\endsubhead
In the setup of 18.2 assume that $G^0$ is semisimple. Let $G_1$ be the subgroup of $G$ 
generated by the connected component that contains $su$. Clearly, $su$ is isolated in 
$G_1$. Applying 18.10 to $G_1,s,u$ instead of $G,s,u$ we see that $s$ is isolated in 
$G_1$. Hence $s$ is isolated in $G$.

\subhead 18.12\endsubhead
In the setup of 18.2 let $\p':G@>>>\bG,\bs,\bu$ be as in 17.12. Then the 
$\bG^0$-conjugacy class of $\bs\bu=\bu\bs$ is isolated in $\bG$. Applying 18.11 to 
$\bG,\bs,\bu$ instead of $G,s,u$ we see that $\bs$ is isolated in $\bG$. Using now 
2.3(a) we see that $s$ is isolated in $G$. Proposition 18.2 is proved.

\head 19. Properties of cuspidal local systems\endhead
\subhead 19.1\endsubhead
Let $s\in G$ be semisimple and let $\fc$ be a unipotent $Z_G(s)^0$-conjugacy class in 
$Z_G(s)$. Assume that the unique $G^0$-conjugacy class $\boc$ that contains $s\fc$ is 
isolated in $G$. Let $G_1=\{g\in Z_{G^0}(s);g\boc g\i=\boc\}$ (a subgroup of 
$Z_{G^0}(s)$ containing $Z_G(s)^0$). Let $\ti{\boc}$ be the variety of orbits for the 
$Z_G(s)^0$-action $z:(y,u)\m(yz\i,zuz\i)$ on $G^0\T\fc$. Then $(y,u)\m ysuy\i$ is a 
finite principal covering $\p:\ti{\boc}@>>>\boc$ with group $G_1/Z_G(s)^0$. Let $\ff$ 
be a cuspidal local system on $\fc$. Let $\tiff$ be the local system on $\ti{\boc}$ 
whose inverse image under $G^0\T\fc@>>>\ti{\boc}$ is $\bbq\bxt\ff$. We show that

(a) $\p_!\tiff$ {\it is a cuspidal local system on $\boc$.}
\nl
Since $\p_!\tiff$ is clearly a $G^0$-equivariant local system, it is enough to show 
that the local system $j^*\p_!\tiff$ on $\fc$ is cuspidal relative to $Z_G(s)$ (here 
$j:\fc@>>>\boc$ is $u\m su$); see 17.3(a). From the definitions we see that 
$j^*\p_!\tiff\cong\op_{g_1}\Ad(g_1)^*\ff$ where $g_1$ runs over a set of 
representatives for the $Z_G(s)^0$-cosets in $G_1$. Clearly, each $\Ad(g_1)^*\ff$ is a 
cuspidal local system on $\fc$ and (a) follows.

\subhead 19.2\endsubhead
{\it Let $u\in G$ be unipotent, quasi-semisimple. Then $Z_{G^0}(u)$ is connected.}
\nl
(See \cite{\DM, 1.28}).

\subhead 19.3\endsubhead
Let $P$ be a parabolic of $G^0$ and let $x\in N_GP,v\in U_P,x'=xv$. We show that

(a) {\it there exists $v'\in U_P$ such that $x'_s=v'x_sv'{}\i$,
$x'_u=v'x_uv'{}\i\mod Z_G(x'_s)^0$.}
\nl
By 1.4(a) we can find Levi subgroups $L.L'$ of $P$ such that $x_s\in N_GL$,
$x'_s\in N_GL'$. Applying the canoical projection $p:N_GP@>>>N_GP/U_P$ to
$x_sx_uv=x'_sx'_u$ we obtain $p(x_s)p(x_u)=p(x'_s)p(x'_u)$. Using the uniqueness of the
Jordan decomposition in $N_GP/U_P$ we get $p(x_s)=p(x'_s)$. We can find $v'\in U_P$ 
such that $L'=v'Lv'{}\i$. Then $v'x_sv'{}\i,x'_s$ are elements of $N_GL'\cap N_GP$ with
the same image under $p$ hence, by 1.26(a), we have $v'x_sv'{}\i=x'_s$. We must show 
that $x_u\i v'{}\i x'_uv'\in Z_G(x_s)^0$. We have $x'_u=v'x_s\i v'{}\i x_sx_uv$ hence 
$x_u\i v'{}\i x'_uv'=x_u\i x_s\i v'{}\i x_sx_uvv'\in U_P$ since $x=x_sx_u\in N_G(U_P)$.
Since $x_u\i v'{}\i x'_uv'\in Z_G(x_s)$ we see that
$x_u\i v'{}\i x'_uv'\in Z_G(x_s)\cap U_P\sub Z_G(x_s)^0$ (we use 1.11). 

\subhead 19.4\endsubhead
Let $C$ be an isolated stratum of $G$ and let $\ce\in\cs(C)$. We show that conditions 
(i),(ii) below are equivalent.

(i) $\ce$ is a cuspidal local system on $C$;

(ii) for any $G^0$-conjugacy class $\boc$ in $C$, $\ce|_\boc$ is a cuspidal local 
system on $\boc$.
\nl
Let $P$ be a parabolic of $G^0$ with $P\ne G^0$ and let $R$ be a $U_P$-coset in $N_GP$.
By 19.3(a), the semisimple part of any element of $R$ is contained in a fixed
$G^0$-conjugacy class. Hence $R$ is contained in a union of finitely many
$G^0$-conjugacy classes. Hence $C\cap R$ is contained in a union of finitely many 
$G^0$-conjugacy classes in $C$; this union is necessarily disjoint (as a variety), by 
the definition of $C$. Thus, $C\cap R=\sqc_{i=1}^n(\boc_i\cap R)$ where $\boc_i$ are 
$G^0$-conjugacy classes in $C$. Let $d$ be the dimension of any $G^0$-conjugacy class 
in $C$ minus the dimension of the $P/U_P$-conjugacy class of $R/U_P$ in $N_GP/U_P$. If
(ii) holds then $H_c^d(C\cap R,\ce)=\op_{i=1}^nH_c^d(\boc_i\cap R,\ce|_{\boc_i})=0$ 
hence
(i) holds. Conversely, assume that (i) holds and $\boc$ is a $G^0$-conjugacy class in 
$C$. We must show that $H_c^d(\boc\cap R,\ce|_\boc)=0$ for $R$ as above. We may assume 
that $\boc\cap R\ne\em$ hence $\boc=\boc_i$ for some $i$. We have 
$0=H_c^d(C\cap R,\ce)=\op_{i=1}^nH_c^d(\boc_i\cap R,\ce|_{\boc_i})$ hence each 
$H_c^d(\boc_i\cap R,\ce|_{\boc_i})$ is $0$. In particular,
$H_c^d(\boc\cap R,\ce|_\boc)=0$, as desired.

\subhead 19.5\endsubhead
Let $C$ be an isolated stratum of $G$. Let $\boc$ be a $G^0$-conjugacy class in $C$. 
Let $\cf$ be a $G^0$-equivariant cuspidal local system on $\boc$. Let $D$ be the 
connected component of $G$ that contains $C$ and let $T={}^D\cz_{G^0}^0$. Let 
$\cl\in\cs(T)$. Define $\p:T\T\boc@>>>C$ by $\p(z,c)=zc$. We show that

(a) $\p_!(\cl\bxt\cf)\in\cs(C)$ {\it is a cuspidal local system.}
\nl
Let $\G$ be the set of all $z\in T$ such that $z\boc=\boc$ (a finite group, see 
1.23(a).) Then $\p$ is a finite principal covering with group $\G$. Hence 
$\p_!(\cl\bxt\cf)$ is a local system on $C$. It is immediate that 
$\p_!(\cl\bxt\cf)\in\cs(C)$. We show that it is cuspidal. Let $P$ be a parabolic of 
$G^0$ with $P\ne G^0$ and let $R$ be a $U_P$-coset in $N_GP$. Let $d$ be $\dim\boc$
minus the dimension of the $P/U_P$-conjugacy class of $R/U_P$ in $N_GP/U_P$. We must 
show that $H_c^d(C\cap R,\p_!(\cl\bxt\cf))=0$ or equivalently that
$H_c^d(\p\i(C\cap R),\cl\bxt\cf)=0$. Now $\p\i(C\cap R)=\{(z,c)\in T\T\boc;zc\in R\}$. 
By 19.3(a), the semisimple part of any element of $R$ is contained in a fixed 
$G^0$-conjugacy class. Hence for $(z,c)\in\p\i(C\cap R)$, $(zc)_s$ is contained in a 
fixed semisimple $G^0$-conjugacy class hence $zc$ is contained in a union of finitely 
many $G^0$-conjugacy classes, hence $z$ can take only finitely many values. Thus there 
exist $z_1,z_2,\do,z_m$ in $T\T\boc$ such that
$$\p\i(C\cap R)=\sqc_{i=1}^m\{(z_i,c);c\in\boc\cap z_i\i R\},$$
$$H_c^d(\p\i(C\cap R),\cl\bxt\cf)
\cong\op_{i=1}^m\cl_{z_i}\ot H_c^d(\boc\cap z_i\i R,\cf)=0.$$
This proves (a).

\subhead 19.6\endsubhead
Let $H$ be a connected algebraic group acting transitively on the variety $X$. Assume 
that we are given $\FF_q$-rational structures on $H,X$ compatible with the action. Let 
$F:H@>>>H$, $F:X@>>>X$ be the Frobenius maps. Let $\Up$ be a set of representatives for
the isomorphism classes of irreducible $H$-equivariant local systems $\cf$ on $X$ such
that $F^*\cf\cong\cf$. For any $\cf\in\Up$ we choose $\ph:F^*\cf@>\si>>\cf$. Then
$\c_{\cf,\ph}:X^F@>>>\bbq$ is a function constant on the orbits of $H^F$, independent 
of the choice of $\ph$, up to a non-zero scalar.

\proclaim{Lemma 19.7} $(\c_{\cf,\ph})_{\cf\in\Up}$ is a $\bbq$-basis of the vector
space of functions $X^F@>>>\bbq$ that are constant on the orbits of $H^F$.
\endproclaim
A special case of this (when $H$ is reductive and $X$ is a unipotent conjugacy class in
$H$) is proved in \cite{\CSV, \S24, p.140}. A similar proof works in the general case.
We can find $x\in X^F$. Let $H_x=\{h\in H;hx=x\}$. Associating to $\cf\in\Up$ the stalk
$\cf_x$ (an irreducible $H_x$-module, by the equivariance of $\cf$, on which $H_x$ acts
through its finite quotient $\G=H_x/H_x^0$) gives a bijection between $\Up$ and a set 
$\Up'$ of representatives for the isomorphism classes of irreducible $\bbq[\G]$-modules
$V$ such that there exists an isomorphism $\io_V:V@>>>V$ with 
$\io_V\g=F(\g)\io_V:V@>>>V$ for all $\g\in\G$. (For $V=\cf_x$ we may take $\io_V$ to be
the isomorphism $\cf_x@>>>\cf_x$ induced by $\ph\i$.) Now $F$ acts naturally on $\G$
and, according to \cite{\SPS, 2.7}, 

(a) {\it $H_x@>>>H^F\bsl X^F$, $z\m H^F-\text{orbit of }hx$ where $h\in H,h\i F(h)=z$
\nl
induces a bijection between the set of $F$-twisted conjugacy classes in $\G$ and the 
set $H^F\bsl X^F$ of $H^F$-orbits on $X^F$.} 
\nl
Via this bijection, giving a function $X^F@>>>\bbq$ that is constant on $H^F$-orbits is
the same as giving a function $\G@>>>\bbq$ that is constant on $F$-twisted conjugacy 
classes in $\G$. If $\cf\in\Up$ and $V=\cf_x$, then the function 
$\c_{\cf,\ph}:X^F@>>>\bbq$ corresponds to the function 
$$\align&\g\m\tr(\cf_{hx}@>F(h)\i>>\cf_x@>\ph>>\cf_x@>h>>\cf_{hx})=
\tr(\cf_x@>F(h)\i h>>\cf_x@>\ph>>\cf_x)\\&=\tr(\io_V\i\g\i,V)\endalign$$
where $h\in H$ is such that $h\i F(h)\in H_x$ has image $\g$ in $\G$. (We use the fact 
that, for any $h'\in H,y\in X$, the compositions 
$$\cf_{F(y)}@>\ph>>\cf_y@>h'>>\cf_{h'y},\qua 
\cf_{F(y)}@>F(h')>>\cf_{F(h')F(y)}@>\ph>>\cf_{h'y}$$
coincide.) It is then enough to show that the functions $\g\m\tr(\io_V\i\g\i,V)$ (for 
various $V\in\Up'$) form a basis for the vector space of functions $\G@>>>\bbq$ that 
are constant on $F$-twisted conjugacy classes. This follows from a variant of the Schur
orthogonality relations. (It also follows from 20.4(f) applied to the group algebra 
$\EE=\bbq[\G]$; in this case all elements of $\G$ are effective, see 20.4.)

\subhead 19.8\endsubhead
In the remainder of this section we assume that $\kk$ is an algebraic closure of a 
finite field $\FF_q$ and that $G$ has a fixed $\FF_q$-rational structure with Frobenius
map $F:G@>>>G$. 

Let $C$ be an isolated stratum of $G$ such that $F(C)=C$. Let $\ci$ be a set of 
representatives for the isomorphism classes of irreducible local systems $\ce$ in 
$\cs(C)$ such that $F^*\ce\cong\ce$; for each $\ce\in\ci$ we choose an isomorphism 
$\ph:F^*\ce@>\si>>\ce$. We show:

(a) {\it the functions $\c_{\ce,\ph}$ where $\ce\in\ci$ form a basis of the vector 
space of functions $C^F@>>>\bbq$ that are constant on the $G^{0F}$-conjugacy classes in
$C^F$.}
\nl
Let $D$ be the connected component of $G$ that contains $C$. Since $\ci$ is a finite 
set, we can find an integer $n\ge 1$, invertible in $\kk$, such that

($*$) $\ce\in\ci\implies \ce\in\cs_n(C)$, (see 5.2),

($**$) $z\in{}^D\cz_{G^0}^0, F(z)=z \implies z^n=1$.
\nl
Applying Lemma 19.7 to the transitive action 5.2(a) (with $n$ as above) of 
$H={}^D\cz_{G^0}^0\T G^0$ on $C$, we see that the functions $\c_{\ce,\ph}$ where $\ce$ 
runs over the elements of $\ci$ that belong to $\cs_n(C)$ (or equivalently, $\ce$ runs 
over $\ci$, see ($*$)) form a basis of the vector space of functions $C^F@>>>\bbq$ that
are constant on the $H^F$-orbits in $C^F$. By ($**$), the $H^F$-orbits on $C^F$ are the
same as the $G^{0F}$-conjugacy classes in $C^F$. This proves (a).

\subhead 19.9\endsubhead
For any isolated stratum $C$ (resp. isolated $G^0$-conjugacy class $\boc$) in $G$ such 
that $F(C)=C$ (resp. $F(\boc)=\boc$) let $\cc_{G^0}(C)$ (resp. $\cc_{G^0}(\boc)$) be 
the subspace of the vector space of functions $C^F@>>>\bbq$ (resp. $\boc^F@>>>\bbq$)
spanned by the functions $\c_{\cf,\e}$ where $\cf$ runs through a set of 
representatives for the isomorphism classes of irreducible cuspidal local systems on 
$C$ (resp. on $\boc$) such that $F^*\cf\cong\cf$ and $\e:F^*\cf@>\si>>\cf$ is a fixed 
isomorphism. Clearly, the subspace $\cc_{G^0}(C)$ (resp. $\cc_{G^0}(\boc)$) is 
independent of choices. From 19.8(a) (resp. Lemma 19.7) we see that the functions 
$\c_{\cf,\e}$ (as above) form a basis of $\cc_{G^0}(C)$ (resp. $\cc_{G^0}(\boc)$): they
are part of a basis of the vector space of all functions $C^F@>>>\bbq$ (resp. 
$\boc^F@>>>\bbq$) that are constant on $G^{0F}$-conjugacy classes. From the definitions
we see that:

(a) {\it if $\cf$ is a (not necessarily irreducible) cuspidal local system on $C$ 
(resp. $\boc$) and $\e:F^*\cf@>>>\cf$ is an isomorphism then $\c_{\cf,\e}$ belongs to 
$\cc_{G^0}(C)$ (resp. $\cc_{G^0}(\boc)$).}

\subhead 19.10\endsubhead
Let $C$ be an isolated stratum of $G$ such that $F(C)=C$. For any $G^0$-conjugacy class
$\boc$ in $C$ such that $F(\boc)=\boc$ then $f\m f|_{\boc^F}$ defines a linear map 

(a) $\cc_{G^0}(C)@>>>\cc_{G^0}(\boc)$,
\nl
(To see that this map is well defined, it is enough to show that, if $\cf$ is a 
cuspidal local system on $C$ and $\e:F^*\cf@>>>\cf$ is an isomorphism then 
$\c_{\cf,\e}|_{\boc^F}\in\cc_{G^0}(\boc)$. This follows from 19.4.) We now take the
direct sum of the maps (a) where $\boc$ runs over the $F$-stable $G^0$-conjugacy 
classes in $C$. We show that 

(b) {\it the resulting linear map $\cc_{G^0}(C)@>>>\op_\boc\cc_{G^0}(\boc)$ is an 
isomorphism.}
\nl
It is obvious that this map is injective. To show that it is surjective it is enough to
verify the following statement: 

{\it for any $F$-stable $G^0$-conjugacy class $\boc$ in $C$, any cuspidal local system 
$\cf$ on $\boc$ and any isomorphism $F^*\cf@>>>\cf$, there exists $f\in\cc_{G^0}(C)$ 
such that $f|_{\boc^F}=\c_{\cf,\e}$ and $f|_{\boc'{}^F}=0$ for any $F$-stable 
$G^0$-conjugacy class $\boc'$ in $C$ with $\boc'\ne\boc$.}
\nl
Let $D$ be the connected component of $G$ that contains $C$. Let $T={}^D\cz_{G^0}^0$.
Let $\cj$ be a set of representatives for the isomorphism classes of local systems 
$\cl$ of rank $1$ in $\cs(T)$ such that $F^*\cl\cong\cl$. For each $\cl\in\cj$ there is
a unique isomorphism $\ph_\cl:F^*\cl@>>>\cl$ which induces the identity map on the 
stalk of $\cl$ at $1$. Then $\th_\cl=\c_{\cl,\ph_\cl}$ is a character $T^F@>>>\bbq^*$ 
and $\cl\m\th_\cl$ is a bijection $\cj@>\si>>\Hom(T^F,\bbq^*)$. Define 
$\p:T\T\boc@>>>C$ by $\p(z,c)=zc$. For any $\cl\in\cj$, $\p_!(\cl\bxt\cf)$ is naturally
isomorphic with its inverse image under $F$ (using $\ph_\cl\bxt\e$); let 
$\c_{\p_!(\cl\bxt\cf),?}:C^F@>>>\bbq$ be the corresponding characteristic function. 
From 19.5(a) we see that $\c_{\p_!(\cl\bxt\cf),?}\in\cc_{G^0}(C)$. From the definitions
we have 

$\c_{\p_!(\cl\bxt\cf),?}(x)=\sum_{z\in T^F,c\in\boc^F;zc=x}\th_\cl(z)\c_{\cf,\e}(c)$
\nl
for $x\in C^F$. Let $f=\sum_{\cl\in\cj}\c_{\p_!(\cl\bxt\cf),?}\in\cc_{G^0}(C)$. For 
$x\in C^F$ we have
$$f(x)=\sum_{z\in T^F,c\in\boc^F;zc=x}\sum_{\cl\in\cj}\th_\cl(z)\c_{\cf,\e}(c)=
\sum_{z\in T^F,c\in\boc^F;zc=x}|T^F|\d_{z,1}\c_{\cf,\e}(c).$$
Thus $f(x)=|T^F|\c_{\cf,\e}(x)$ if $x\in\boc^F$ and $f(x)=0$ if $x\in C^F-\boc^F$.
This completes the proof of (b).

\subhead 19.11\endsubhead
{\it If $E$ is the set of unipotent quasi-semisimple elements in some $F$-stable 
connected component of $G$ that contains unipotent elements then $E^F$ is a single 
$G^{0F}$-conjugacy class.}
\nl
This follows from the fact that $E$ is a homogeneous $G^0$-space (see 1.9(a)) defined
over $\FF_q$ in which the isotropy group of any point is connected (see 19.2).

\subhead 19.12\endsubhead
Let $\boc$ be a cuspidal $G^0$-conjugacy class in $G$. Let $\bocd=\s_G(x)$ for any 
$x\in\boc$; this is the set of all quasi-semisimple elements $g\in G$ such that for 
some $x\in\boc$ we have $x_s=g_s$, $x_u\in Z_G(g_s)^0g_u$. Let $Z$ be the set of all 
pairs $(s,c)$ where $s\in G$ is semisimple and $c$ is a connected component of $Z_G(s)$
such that there exists a unipotent element $u\in Z_G(s)$ with $su\in\boc,u\in c$. We 
have a diagram 
$$\boc@>a>>Z@<b<<\bocd$$
where $a(x)=(x_s,Z_G(x_s)^0x_u),b(g)=(g_s,Z_G(g_s)^0g_u)$. Now $G^0$ acts transitively 
on $\boc,Z,\bocd$ compatibly with $a,b$. For any $y\in G$ let 
$$H_{G^0}(y)=\{h\in G^0;hy_sh\i=y_s,hy_uh\i\in Z_G(y_s)^0y_u\},$$
a closed subgroup of $Z_{G^0}(y_s)$ containing $Z_G(y_s)^0$.

(a) {\it Let $x\in\boc,g\in\bocd$ be such that $a(x)=b(g)$. Let
$H=H_{G^0}(x)=H_{G^0}(g)$ (the stabilizer of $a(x)=b(g)$ in $G^0$). The map 
$Z_{G^0}(x)/Z_{G^0}(x)^0@>>>H/H^0$ induced by $a$ is surjective and the map
$Z_{G^0}(g)/Z_{G^0}(g)^0@>>>H/H^0$ induced by $b$ is an isomorphism.}
\nl
We have $a(x)=b(g)=(x_s,Z_G(x_s)^0x_u)=(g_s,Z_G(g_s)^0g_u)$. Let $h\in H$. Then 
$x,hxh\i$ are elements of $\boc$ with the same semisimple part $x_s$ and their 
unipotent parts are contained in the same connected component of $Z_G(x_s)$. By 17.13, 
there exists $z\in Z_G(x_s)^0$ such that $hx_uh\i=zx_uz\i$. Thus, $h=zh_1$ where 
$h_1\in Z_{G^0}(x)$. We see that $H=Z_G(x_s)^0Z_{G^0}(x)$. This proves the first 
assertion of (a).

If $h\in H$ then $g_u,hg_uh\i$ are unipotent quasi-semisimple elements of $Z_G(g_s)$ 
contained in the same connected component of $Z_G(g_s)$ hence, by 17.13, there exists 
$z\in Z_G(g_s)^0$ such that $hg_uh\i=zg_uz\i$. Thus, $h=zh_1$ where 
$h_1\in Z_{G^0}(g)$. We see that $H=Z_G(g_s)^0Z_{G^0}(g)$. This shows that $b$ induces 
a surjective map $Z_{G^0}(g)/Z_{G^0}(g)^0@>\si>>H/H^0$ and that the group of components
of $H$ is the same as the group of components of 
$Z_{G^0}(g)/(Z_{G^0}(g)\cap Z_G(g_s)^0)=Z_{G^0}(g)/(Z_{Z_G(g_s)^0}(g_u)$ which (by the 
connectedness of $Z_{Z_G(g_s)^0}(g_u)$, see 19.2) is the same as the group of 
components of $Z_{G^0}(g)$. Thus, the surjective map
$Z_{G^0}(g)/Z_{G^0}(g)^0@>\si>>H/H^0$ must be an isomorphism. This proves (a).

\subhead 19.13\endsubhead
In the setup of 19.12, let $\cl$ be an irreducible $G^0$-equivariant local system on 
$Z$. Since $H/H^0=Z_{G^0}(g)/Z_{G^0}(g)^0$ (see 19.12(a)) is commutative (see 17.4(a)),
$\cl$ has rank $1$. Let $\tcl=a^*\cl$, a local system of rank $1$ on $\boc$. We show:

(a) {\it if $\cf$ is a cuspidal local system on $\boc$ then $\cf\ot\tcl$ is a cuspidal
local system  on $\boc$.}
\nl
Let $P$ be a parabolic of $G^0$ with $P\ne G^0$ and let $x\in\boc\cap N_GP$. Let $d$ be
$\dim\boc$ minus the dimension of the $P/U_P$-conjugacy class of $xU_P$ in $N_GP/U_P$. 
We must show that $H_c^d(\boc\cap xU_P,\cf\ot\tcl)=0$. By our assumption we have 
$H_c^d(\boc\cap xU_P,\cf)=0$. Hence it is enough to show that
$\tcl|_{\boc\cap xU_P}\cong\bbq$. Since $\tcl=a^*\cl$, it is enough to show that there 
exists a subvariety $V$ of $Z$ such that 

(b) $a(\boc\cap xU_P)\sub V$, $\cl|_V\cong\bbq$.
\nl
Let $V$ be the $U_P$-orbit of $a(x)$ in $Z$ (for the restriction of the $G^0$-action to
$U_P$). For this $V$ the first assertion of (b) holds by 19.3(a). We now show that for 
this $V$, the second assertion of (b) holds. It is enough to note that $\cl|_V$ is a 
$U_P$-equivariant local system of rank $1$ on the homogeneous $U_P$-space $V$ in which 
the isotropy group of $a(x)$ that is, $U_P\cap Z_G(x_s)$, is connected (we use that 
$x_s$ normalizes $U_P$, see 1.11). Thus (b), hence also (a), are proved. 

\subhead 19.14\endsubhead
We now assume that $F(\boc)=\boc$. Then $Z$ and $\boc$ are defined over $\FF_q$ and we
denote again by $F$ the corresponding Frobenius maps. 

(a) {\it The map $a_0:Z^F@>>>\boc^F$ (restriction of $a:Z@>>>\boc$) is surjective; the 
map $b_0:Z^F@>>>\bocd^F$ (restriction of $b:Z@>>>\bocd$) induces a bijection on the 
sets of $G^{0F}$-orbits.}
\nl
This follows immediately from 19.12(a) and 19.7(a). We have a partition 
$\bocd^F=\sqc\g$
where $\g$ runs over the $G^{0F}$-orbits on $\bocd^F$. For any $\g$ we set 
$Z^F_\g=b_0\i(\g),\boc^F_\g=a_0(Z^F_\g)$. From (a) we see that $Z^F=\sqc_\g Z^F_\g$ is 
the partition of $Z^F$ into $G^{0F}$-orbits and that $\boc^F=\sqc_\g\boc^F_\g$ is a 
partition of $\boc^F$ into non-empty $G^{0F}$-stable subsets.

Let $f\in\cc_{G^0}(\boc)$ and let $\g$ be a $G^{0F}$-orbit on $\bocd^F$. Define 
$f_\g:\boc^F@>>>\bbq$ by $f_\g(x)=1$ if $\g\in\boc^F_\g,f_\g(x)=0$ if 
$\g\in\boc^F-\boc^F_\g$. We show that
$$f_\g f\in\cc_{G^0}(\boc).\tag b$$
We may assume that $f=\c_{\cf,\e}$ where $\cf$ is a cuspidal local system on $\boc$ and
$\e:F^*\cf@>\si>>\cf$ is an isomorphism. Define $\tf_\g:Z^F@>>>\bbq$ by $\tf_\g(z)=1$
if $z\in Z^F_\g$, $\tf_\g(z)=0$ if $z\in Z^F-Z^F_\g$. By 19.7 applied to the 
homogeneous $G^0$-space $Z$, there exist irreducible $G^0$-equivariant local systems 
$\cl^i$, $(i\in[1,m]$ on $Z$ and isomorphisms $e^i:F^*\cl^i@>\si>>\cl^i$ such that
$\tf_\g=\sum_{i\in[1,m]}c_i\c_{\cl^i,e^i}$ where $c_i\in\bbq$. Each $\cl^i$ has 
rank $1$. Composing with $a_0:\boc^F@>>>Z^F$ we obtain 
$\tf_\g\circ a_0=\sum_{i\in[1,m]}c_i\c_{\cl^i,e^i}\circ a_0$, that is
$f_\g=\sum_{i\in[1,m]}c_i\c_{\tcl^i,\ti e^i}$ where $\tcl^i=a^*\cl^i$ and
$\ti e^i:F^*\tcl^i@>\si>>\tcl^i$ is induced by $e^i$. Hence 
$$f_\g f=\sum_{i\in[1,m]}c_i\c_{\tcl^i,\ti e^i}\c_{\cf,\e}=
\sum_{i\in[1,m]}c_i\c_{\cf\ot\tcl^i,\e\ot\ti e^i}.$$
Using 19.13(a), we see that this belongs to $\cc_{G^0}(\boc)$. This proves (b).

For any $G^{0F}$-orbit $\g$ on $\bocd^F$ we set
$$\cc_{G^0,\g}(\boc)=\{f\in\cc_{G^0}(\boc);f=0\text{ on }\boc^F-\boc^F_\g\}.$$
From (b) we see that:
$$\cc_{G^0}(\boc)=\op_\g\cc_{G^0,\g}(\boc).\tag c$$

\subhead 19.15\endsubhead
We now fix $g\in\bocd^F$. We set $s=g_s$. Let

$\fc=\{u\in Z_G(s)^0g_u;u\text{ unipotent, }su\in\boc\}$. 
\nl
From the definitions we see 
that $\fc\ne\em$ and from 17.13 we see that $\fc$ is a single 
(unipotent) $Z_G(s)^0$-conjugacy class in $Z_G(s)$. We have $F(\fc)=\fc$ and $\fc$ 
carries some non-zero cuspidal local system (see 17.3(a)). Define $\cc_{Z_G(s)^0}(\fc)$
in terms of $Z_G(s),\fc$ in the same way as $\cc_{G^0}(\boc)$ was defined in terms of
$G,\boc$. For any $f\in\cc_{G^0}(\boc)$ we define $\baf:\fc^F@>>>\bbq$ by 
$\baf(u)=f(su)$. We claim that 

$\baf\in\cc_{Z_G(s)^0}(\fc)$. 
\nl
We may assume that $f=\c_{\cf,\e}$ where $\cf$ is a cuspidal local system on $\boc$ and
$\e:F^*\cf@>\si>>\cf$ is an isomorphism. Let $\cf'=j^*\cf$ where $j:\fc@>>>\boc,u\m su$
and let $\e':F^*\cf'@>\si>>\cf'$ be the isomorphism induced by $\e$. By 17.3(a), $\cf'$
is a cuspidal local system on $\fc$ and by 19.9(a) applied to $\fc$ instead of $\boc$ 
we see that $\c_{\cf',\e'}\in\cc_{Z_G(s)^0}(\fc)$. Clearly, $\c_{\cf',\e'}=\baf$ and 
our claim is verified.

From 17.13 it follows that $\fc^F$ is stable under conjugation by 
$H_{G^0}(g)^F$. It is clear that for $f,\baf$ as above, $\baf$ is constant on any 
$H_{G^0}(g)^F$-conjugacy class in $\fc^F$. Thus we have a well defined linear map
$$\cc_{G^0}(\boc)@>>>\cc_{H_{G^0}(g)}(\fc),f\m\baf\tag a$$ 
where $\cc_{H_{G^0}(g)}(\fc)$ is the space of functions in $\cc_{Z_G(s)^0}(\fc)$ that 
are constant on any $H_{G^0}(g)^F$-conjugacy class in $\fc^F$. 

We show that the map (a) is surjective. Now $\cc_{H_{G^0}(g)}(\fc)$ is spanned by 
functions $f':\fc^F@>>>\bbq$ of the form
$$f'(u)=|Z_G(s)^{0F}|\i\sum_{y\in H_{G^0}(g)^F}\c_{\ff,\e}(y\i uy)$$
where $\ff$ is a cuspidal local system on $\fc$ and $\e:F^*\ff@>\si>>\ff$ is an 
isomorphism. It is enough to show that any such $f'$ is in the image of the map (a). 
Let $\p:\ti{\boc}@>>>\boc,\tiff$ be as in 19.1. Now $\e$ induces an isomorphism 
$\ti\e:F^*\p_!\tiff@>>>\p_!\tiff$. Let $f=\c_{\p_!\tiff,\ti\e}:\boc^F@>>>\bbq$. Since 
$\p_!\tiff$ is a cuspidal local system on $\boc$ (see 19.1(a)) we see that 
$f\in\cc_{G^0}(\boc)$. From the definitions we have 
$$\baf(u)=f(su)=|Z_G(s)^{0F}|\i\sum_{y\in G^{0F},u'\in\fc^F;ysu'y\i=su}\c_{\ff,\e}(u')
$$
for $u\in\fc^F$. For each $y,u'$ in the sum we have $y\in Z_{G^0}(s)^F$ and
$yu'y\i=u$ hence $y\in H_{G^0}(g)^F, u'=y\i uy$; hence
$$\baf(u)=|Z_G(s)^{0F}|\i\sum_{y\in H_{G^0}(g)^F}\c_{\ff,\e}(y\i uy).$$
We see that $\baf=f'$. Thus the surjectivity of the map (a) is established.

Next we note that, if $\g'$ is a $G^{0F}$-orbit on $\bocd^F$ that does not contain $g$,
then the restriction of the map (a) to $\cc_{G^0,\g'}(\boc)$ is $0$. Using this and the
direct sum decomposition 19.14(c) we deduce that, if $\g$ is the $G^{0F}$-orbit on
$\bocd^F$ that contains $g$, then (a) restricts to a surjective linear map
$$\cc_{G^0,\g}(\boc)@>>>\cc_{H_{G^0}(g)}(\fc).\tag b$$ 
We show that this map is injective. Let $f,f'\in\cc_{G^0,\g}(\boc)$ be such that 
$f(su)=f'(su)$ for any $u\in\fc^F$. Since $\boc^F_\g=\{ysuy\i;y\in G^{0F},u\in\fc^F\}$ 
and $f,f'$ are constant on $G^{0F}$-conjugacy classes it follows that $f=f'$ on 
$\boc^F_\g$. Since $f,f'$ are $0$ on $\boc^F-\boc^F_\g$ it follows that $f=f'$, as 
desired. We see that

(c) {\it the map (b) is an isomorphism.}

\head 20. Twisted group algebras\endhead
\subhead 20.1\endsubhead
Let $\G$ be a finite group. Let $\EE$ be a finite dimensional $\bbq$-vector space with 
a direct sum decomposition $\EE=\op_{w\in\G}\EE_w$ with $\dim\EE_w=1$ for all $w$. 
Assume that on $\EE$ we are given an associative algebra structure with $1$ such that 
$\EE_w\EE_y=\EE_{wy}$ for any $w,y\in\G$. Then $1\in\EE_1-\{0\}$. We choose a basis 
$\{b_w;w\in\G\}$ of $\EE$ such that $b_w\in\EE_w$ for all $w$. Each $b_w$ is 
invertible. We have $b_wb_y=\l(w,y)b_{wy},b_y\i b_w\i=\l(w,y)\i b_{wy}\i$ with
$\l(w,y)\in\bbq^*$ for any $w,y\in\G$. We show that

(a) {\it the algebra $\EE$ is semisimple.} 
\nl
Let $M$ be a finite dimensional $\EE$-module and let $M'$ be an $\EE$-submodule of $M$.
We must show that there exists an $\EE$-submodule of $M$ complementary to $M'$. Let 
$\p:M@>>>M'$ be a $\bbq$-linear map such that $\p(m')=m'$ for all $m'\in M'$. Define a 
$\bbq$-linear map $\ti\p:M@>>>M'$ by $\ti\p(m)=|\G|\i\sum_{w\in\G}b_w\i\p(b_wm)$. For 
$m'\in M'$ we have $\ti\p(m')=|\G|\i\sum_{w\in\G}b_w\i b_wm'=m'$. We show that $\ti\p$ 
is $\EE$-linear. It is enough to show that $b_y\i\ti\p(b_ym)=\ti\p(m)$ for $m\in M$,
$y\in\G$. We have 
$$\align&|\G|b_y\i\ti\p(b_ym)=b_y\i\sum_wb_w\i\p(b_wb_ym)=\sum_wb_y\i b_w\i 
\p(\l(w,y)b_{wy}m)\\&=\sum_w\l(w,y)\i b_{wy}\i\p(\l(w,y)b_{wy}m)
=\sum_wb_{wy}\i\p(b_{wy}m)=|\G|\ti\p(m),\endalign$$
as desired. Now $\Ker\ti\p$ is an $\EE$-submodule of $M$ complementary to $M'$. This 
proves (a).

\subhead 20.2\endsubhead
Let $V,V'$ be two simple $\EE$-modules. Let $t:V@>>>V,t':V'@>>>V'$ be $\bbq$-linear 
invertible maps. Let $N=|\G|\i\sum_{w\in\G}\tr(b_wt,V)\tr(t'{}\i b_w\i,V')$. We show:

(a){\it $N=0$ if $V,V'$ are not isomorphic $\EE$-modules and $N=1$ if $V=V',t=t'$.}
\nl
In the case where $\l(w,y)=1$ for all $w,y$ this is just Schur's orthogonality formula.
The proof in the general case is similar. Let $(e_i)_{i\in I}$ be a basis of $V$ and 
let $(e'_h)_{h\in I'}$ be a basis of $V'$. For $w\in\G$ we define 
$\a_{ij}^w,\b_{hk}^w\in\bbq$ by $b_w(e_i)=\sum_{j\in I}\a_{ij}^we_j$, 
$b_w\i(e'_h)=\sum_{k\in I'}\b_{hk}^we'_k$. Define $\xi_{ij},\z_{hk}\in\bbq$ by 
$t(e_i)=\sum_{j\in I}\xi_{ij}e_j$, $t'{}\i(e'_h)=\sum_{k\in I'}\z_{hk}e'_k$. We have
$$N=|\G|\i\sum_{w\in\G}\sum_{i,j,k,h}\xi_{ij}\a_{ji}^w\z_{hk}\b_{kh}^w.\tag b$$
For a $\bbq$-linear map $f:V@>>>V'$ we define $\tf:V@>>>V'$ by 
$\tf(v)=\sum_{w\in\G}b_w\i f(b_wv)$. As in the proof of 20.1(a) we see that $\tf$ is 
$\EE$-linear. For $i\in I,h\in I'$ define a linear map $f:V@>>>V'$ by 
$f(e_j)=\d_{ij}e'_h$ for all $j$. We have
$$\align&\tf(e_u)=\sum_wb_w\i f(b_we_u)=\sum_wb_w\i f(\sum_j\a_{uj}^we_j)=
\sum_wb_w\i\sum_j\a_{uj}^w\d_{ij}e'_h\\&=
\sum_w\a_{ui}^wb_w\i e'_h=\sum_{w,k}\a_{ui}^w\b_{hk}^we'_k.\endalign$$
If $V,V'$ are not isomorphic $\EE$-modules then, by Schur's lemma, we have $\tf=0$ 
hence $\sum_w\a_{ui}^w\b_{hk}^w=0$ for any $u,i,h,k$ and $N=0$ as desired. Assume now 
that $V=V',t=t'$. We may assume that $I=I'$, $e_i=e'_i$. If $f,\tf$ are as above then, 
by Schur's lemma, $\tf$ is a multiple of $1$. Since 
$\tr(\tf,V)=|\G|\tr(f,V)=|\G|\tr(f,V)=\d_{ih}|\G|$ we have $\tf=\d_{ih}|\G|n\i 1$ where
$n=\dim V$. Hence $\sum_w\a_{ui}^w\b_{hk}^w=\d_{ih}\d_{uk}|\G|n\i$ for any $u,i,h,k$ 
and 
$$N=\sum_{i,j,k,h}\xi_{ij}\d_{ik}\d_{jh}\z_{hk}n\i=
\sum_{i,j}\xi_{ij}\z_{ji}n\i=\sum_in\i=1,$$
as desired.

\subhead 20.3\endsubhead
Since $\EE$ is semisimple, we have an algebra isomorphism 

(a) $\EE@>>>\op_{i=1}^{r'}\End(V_i)$,
\nl
$e\m(f_i^e)$ where $f_i^e:V_i@>>>V_i$ takes $v$ to $ev$; here $V_i,(i\in[1,r'])$ is a 
set of representatives for the isomorphism classes of simple $\EE$-modules.

Let $\io:\EE@>>>\EE$ be an automorphism of the algebra $\EE$. We may assume that for 
$i\in[1,r]$ the following property holds: 

($*$) {\it there exists a $\bbq$-linear isomorphism $\io_i:V_i@>>>V_i$ such that 
$\io_i(ev)=\io(e)\io_i(v)$ for all $e\in\EE,v\in V$}
\nl
and that for $i>r$ this property does not hold. We choose for each $i\in[1,r]$ an
isomorphism $\io_i:V_i@>>>V_i$ as above (it is unique up to a non-zero scalar). We show
that, for any $w,w'$ in $\G$,

(b) $\sum_{i=1}^r\tr(b_w\io_i,V_i)\tr(\io_i\i b_{w'}\i,V_i)$ {\it is equal to the trace
of the linear map $\k:\EE@>>>\EE,e\m b_{w'}\i\io\i(e)b_w$.}
\nl
Let $\t:\op_{i=1}^{r'}\End(V_i)@>>>\op_{i=1}^{r'}\End(V_i)$ be the linear map which 
corresponds to $\k:\EE@>>>\EE$ under the isomorphism (a). For $i\in[1,r]$, $\t$ 
restricts to a linear map $\t_i:\End(V_i)@>>>\End(V_i)$, while for $i>r$, $\t$ maps the
summand $\End(V_i)$ to a different summand. Hence 
$\tr(\k,\EE)=\sum_{i=1}^r\tr(\t_i,\End(V_i))$. For $i\in[1,k]$, $\t_i$ takes 
$f\in\End(V_i)$ to $v\m b_{w'}\i\io_i\i(f(\io_i(b_wv)))$ hence 

$\tr(\t_i,\End(V_i)=\tr(b_{w'}\i\io_i\i,V_i)\tr(\io_i b_w,V_i)$;
\nl
(b) follows.

\subhead 20.4\endsubhead
We now assume that $\io:\EE@>>>\EE$ in 20.3 satisfies $\io(\EE_w)=\EE_{F(w)}$ for all 
$w$, where $F:\G@>>>\G$ is a group isomorphism. 
For $x\in\G$, let $\G_x=\{y\in\G;F\i(y)xy\i=x\}$; we define 
$\g_x:\G_x@>>>\bbq^*$ by $\io\i(b_y)b_x=\g_x(y)b_xb_y$. We show that $\g_x$ is a group 
homomorphism. Let $z,z'\in\G_x$. We have $b_zb_{z'}=ub_{zz'}$ with $u\in\bbq^*$. We 
have
$$\align&\io\i(b_{zz'})b_x=u\i\io\i(b_zb_{z'})b_x=u\i\io\i(b_z)\g_x(z')b_xb_{z'}\\&=
u\i\g_x(z')\io\i(b_z)b_xb_{z'}=u\i\g_x(z')\g_x(z)b_xb_zb_{z'}=\g_x(z')\g_x(z)
b_xb_{zz'}\endalign$$
hence $\g_x(zz')=\g_x(z')\g_x(z)$, as desired.

An element $x\in\G$ is said to be {\it effective} if $\g_x$ is identically $1$. For
$x,y\in\G,z\in\G_x$ we have $yzy\i\in\G_{F\i(y)xy\i}$ and 
$\g_x(z)=\g_{F\i(y)xy\i}(yzy\i)$. It follows that the set of effective elements in $\G$
is a union of $F$-twisted conjugacy classes. We say that an $F$-twisted conjugacy class
in $\G$ is effective if some/any element of it is effective.

(a) {\it If an $F$-twisted conjugacy class C is not effective then for $i\in[1,r]$ and 
$x\in C$ we have} $\tr(\io_ib_x,V_i)=0$.
\nl
Indeed, we can find $y\in\G_x$ such that $\g_x(y)\ne 1$. We have
$$\align&\tr(\io_ib_x,V_i)=\tr(b_y\i\io_ib_xb_y,V_i)
=\g_x(y)\i\tr(b_y\i\io_i\io\i(b_y)b_x,V_i)\\&=\g_x(y)\i\tr(\io_ib_x,V_i).\endalign$$
Thus $(1-\g_x(y)\i)\tr(\io_ib_x,V_i)=0$ and $\tr(\io_ib_x,V_i)=0$ as claimed.

(b) {\it If $i,j\in[1,r]$ then $x\m\tr(b_x\io_i,V_i)\tr(\io_j\i b_x\i,V_j)$ is constant
on any $F$-twisted conjugacy class}.
\nl
Indeed let $y\in\G$. We have $b(F\i(y)xy\i)=c\io\i(b_y)b_xb_y\i$ for some $c\in\bbq^*$.
Hence 
$$\align&\tr(b_{F\i(y)xy\i}\io_i,V_i)\tr(\io_j\i b_{F\i(y)xy\i}\i,V_j)\\&=
\tr(c\io\i(b_y)b_xb_y\i\io_i,V_i)\tr(\io_j\i c\i b_yb_x\i\io\i(b_y)\i,V_j)\\&
=\tr(\io\i(b_y)b_x\io_i\io\i(b_y)\i,V_i)\tr(\io\i(b_y)\io_j\i b_x\i\io\i(b_y)\i,V_j)\\&
=\tr(b_x\io_i,V_i)\tr(\io_j\i b_x\i,V_j)\endalign$$
and (b) follows.

Let $\baG$ be a set of representatives for the effective $F$-twisted conjugacy classes
in $\G$. We rewrite 20.2(a) for $V=V_i,V'=V_j,t=\io_i,t'=\io_j$ where $i,j\in[1,r]$, 
taking into account (a),(b):
$$\sum_{x\in\baG}|\G_x|\i\tr(b_x\io_i,V_i)\tr(\io_j\i b_x\i,V_j)=\d_{ij}.\tag c$$
We rewrite 20.3(b) for $w,w'\in\baG$ as follows
$$\sum_{i=1}^r\tr(b_w\io_i,V_i)\tr(\io_i\i b_{w'}\i,V_i)=\d_{w,w'}|\G_w|.\tag d$$
Indeed, it is enough to show that the trace of $\k$ in 20.3(b) is in our case 
$\d_{w,w'}|\G_w|$. Now that trace is $\sum_{y\in\G;w'{}\i F\i(y)w=y}c_y$ where
$c_y\in\bbq^*$ is defined by \lb $b_{w'}\i\io\i(b_y)b_w=c_yb_y$. If $w'\ne w$, the last
sum is empty so its value is $0$. If $w'=w$, the last sum is $\sum_{y\in\G_w}\g_w(y)$ 
and this equals $|\G_w|$ since $\g_w$ is identically $1$.

(e) {\it The matrix $(\tr(b_x\io_i,V_i))_{i\in[1,r],x\in\baG}$ is square and 
invertible.}
\nl
Indeed, from (c),(d) we see that this matrix has a left inverse and a right inverse. 
The same argument shows that 

(f) {\it the matrix $(\tr(\io_i\i b_x\i,V_i))_{i\in[1,r],x\in\baG}$ is square and 
invertible.}
\nl
In particular,

(g) $|\baG|=r$.

\head 21. Bases\endhead
\subhead 21.1\endsubhead  
If $L$ is a Levi of a parabolic of $G^0$, let $N_G^\bul L$ be the set of all 
$g\in N_GL$ such that for some parabolic $P$ of $G^0$ with Levi $L$ we have 
$g\in N_GP$. Then $N_G^\bul L$ is a union of connected components of $N_GL$.

Let $\fA_G$ be the set of all pairs $(L,\fc)$ where $L$ is a Levi subgroup of some 
parabolic of $G^0$ and $\fc$ is a unipotent cuspidal $L$-conjugacy class in 
$N_G^\bul L$.

Let $G_{un}$ be the set of unipotent elements in $G$. 

\subhead 21.2\endsubhead
Let $L$ be a Levi of a parabolic $P$ of $G^0$ and let $g\in\tL=N_GL\cap N_GP$. Then

(a) {\it $g$ is quasi-semisimple in $G$ if and only if $g$ is quasi-semisimple in}
$\tL$. 
\nl
See \cite{\DM, 1.10}.

\subhead 21.3\endsubhead
For a fixed $g\in G$, let $\fR$ be the set of all $\uL$ such that $\uL$ is a Levi of a 
parabolic of $G^0$, $g\in N_G^\bul\uL$ and $g$ is isolated in $N_G\uL$; let $\fR'$ be 
the set of all $L$ such that $L$ is a Levi of a parabolic of $Z_G(g_s)^0$ and 
$g\in N_{Z_G(g_s)}^\bul L$. We show:

(a) $\uL\m a(\uL)=\uL\cap Z_G(g_s)^0,L\m b(L)=Z_{G^0}((\cz_L^0\cap Z_G(g))^0)$ 
{\it define inverse bijections $\fR\lra\fR'$.}
\nl
Let $L\in\fR'$. Set $\uL=b(L)$. Then $\uL$ is a Levi of a parabolic of $G^0$. Clearly, 
$g\in N_G\uL$. If $\c:\kk^*@>>>(\cz_L^0\cap Z_G(g))^0$ is general enough, then

$Z_{Z_G(g_s)^0}(\c(\kk^*))=Z_{Z_G(g_s)^0}(\cz_L^0\cap Z_G(g))^0=L$ (see 1.10) and
$Z_{G^0}(\c(\kk^*))=Z_{G^0}(\cz_L^0\cap Z_G(g))^0=\uL$.
\nl
Then $\uL$ is a Levi of the parabolic $P=P_\c$ (see 1.16). We have $g\c(t)g\i=\c(t)$ 
for all $t\in\kk^*$, hence $gPg\i=P$. Thus $g\in N_G^\bul\uL$. We have
$$\uL\cap Z_G(g_s)^0=Z_{G^0}(\c(\kk^*))\cap Z_G(g_s)^0=Z_{Z_G(g_s)^0}(\c(\kk^*))=L.$$
From $L=Z_G(g_s)^0\cap\uL$ we see that $L=Z_{\uL}(g_s)^0$. Hence 
$\cz_{Z_{\uL}(g_s)^0}^0=\cz_L^0$,
$$T_{N_G\uL}(g)=(\cz_{Z_{\uL}(g_s)^0}^0\cap Z_G(g))^0
=(\cz_L^0\cap Z_G(g))^0\sub\cz_{\uL}^0.$$
By 2.2(ii), we see that $g$ is isolated in $N_G\uL$. Hence $\uL\in\fR$ and $b$ is well 
defined.

Conversely, let $\uL\in\fR$. Set $L=a(\uL)$. Let $P$ be a parabolic of $G^0$ with Levi 
$\uL$ such that $g\in N_G\uL\cap N_GP$. Let $Q=P\cap Z_G(g_s)^0$. By 1.12, $Q$ is a 
parabolic of $Z_G(g_s)^0$ with Levi $L$. Clearly, $g\in Z_G(g_s)$ normalizes $L$ and 
$Q$. Hence $L\in\fR'$ and $a$ is well defined. Since $g$ is isolated in $N_G\uL$, we 
have (see 2.2(iii)):
$$(\cz_{\uL}^0\cap Z_G(g))^0=T_{N_G\uL}(g)=(\cz_{Z_{\uL}(g_s)^0}^0\cap Z_G(g))^0=
(\cz_L^0\cap Z_G(g))^0$$
hence (using 1.10): 
$$\uL=Z_{G^0}((\cz_{\uL}^0\cap Z_G(g))^0)=Z_{G^0}((\cz_L^0\cap Z_G(g))^0)=b(L)=ba(\uL).
$$
Thus, $ba=1$. As we have seen above, for $L\in\fR'$ we have $ab(L)=L$. This proves (a).

\subhead 21.4\endsubhead
In the remainder of this section we assume that $\kk$ is an algebraic closure of a 
finite field $\FF_q$ and that $G$ has a fixed $\FF_q$-rational structure with Frobenius
map $F:G@>>>G$. 

Let $\car$ be the set of all triples $(g,\uL,\boc)$ where $g\in G^F$ is 
quasi-semisimple, $\uL$ is a Levi of a parabolic of $G^0$, $F(\uL)=\uL$, $\boc$ is a
a cuspidal $\uL$-conjugacy class in $N_G\uL$ such that $F(\boc)=\boc$, 
$\boc\sub N_G^\bul\uL$ and $g\in\s_{N_G\uL}(\boc)$. 

Let $\car'$ be the set of all triples $(g,L,\fc)$ where $g\in G^F$ is quasi-semisimple,
$L$ is a Levi of a parabolic of $Z_G(g_s)^0$, $F(L)=L$, and $\fc$ is a unipotent 
cuspidal $L$-conjugacy class in $N_{Z_G(g_s)}L$ such that $F(\fc)=\fc$, $\fc\sub g_uL$,
$\fc\sub N_{Z_G(g_s)}^\bul L$ (we have automatically $g\in N_{Z_G(g_s)}L$). Define 
$\ti a:\car@>>>\car'$ by

(a) $(g,\uL,\boc)\m(g,L,\fc),L=a(\uL),
\fc=\{u\in g_uL;u\text{ unipotent,}g_su\in\boc\}$.
\nl
To see that $(g,L,\fc)\in\car'$ we note that $\fc$ is a single $L$-conjugacy class, by
17.13 applied to $N_G\uL$ instead of $G$; also from 20.7(a) we have
$g\in N_{Z_G(g_s)}^\bul L,\fc\sub g_uL$ hence $\fc\sub N_{Z_G(g_s)}^\bul L$.

Define $\ti b:\car'@>>>\car$ by 

(b) $(g,L,\fc)\m(g,\uL,\boc),\uL=b(L),\boc=\cup_{l\in\uL}lg_s\fc l\i$. 
\nl
To see that $(g,\uL,\boc)\in\car$ we note that there exists $x\in\boc$ with 
$x_s=g_s,x_u\in Z_{\uL}(g_s)^0g_u$ and that $g$ is isolated in $N_G\uL$, 
$g\in N_G^\bul\uL$ (see 20.7(a)); it follows that $g\in\s_{N_G\uL}(\boc)$ and that $\boc$ 
is isolated in $N_G\uL$, $\boc\in N_G^\bul\uL$ (we apply Lemma 2.5 with $N_G\uL$ 
instead of $G$ to $g$ and $x$ as above).

From the definition we see that

(c) {\it the maps (a),(b) are inverse bijections $\car\lra\car'$.}

\subhead 21.5\endsubhead
Assume that $(g,\uL,\boc)\in\car,(g,L,\fc)\in\car$ correspond to each other under the 
bijections 20.8(c). Let $\la g\ra=\cup_{l\in\uL^F}lgl\i$. Let 

$\cg=\{y\in Z_{G^0}(g_s)^F;yLy\i=L,y\fc y\i=\fc\}$,

$\cg'=\{y\in G^{0F},y\uL y\i=\uL,y\boc y\i=\boc,y\la g\ra y\i=\la g\ra\}$.
\nl
We show:

(a) $\cg'=\uL^F\cg$.
\nl
Let $y\in\cg$. Since $\fc\sub g_uL$ and $yLy\i=L,y\fc y\i=\fc$, we have 
$\fc\sub yg_uy\i L$. Hence $yg_uy\i,g_u$ are two $F$-stable, unipotent quasi-semisimple
elements of $N_GL$ in the same connected component of $N_GL$. Using 19.11 we have 
$yg_uy\i=lg_ul\i$ for some $l\in L^F$. We have 
$\uL=Z_{G^0}((\cz_L^0\cap Z_G(g))^0)=Z_{G^0}((\cz_L^0\cap Z_G(g_u))^0)$. Hence 
$$\align&y\uL y\i=Z_{G^0}((\cz_{yLy\i}^0\cap Z_G(yg_uy\i))^0)\\&=
Z_{G^0}((\cz_L^0\cap Z_G(lg_ul\i))^0)=Z_{G^0}((\cz_L^0\cap Z_G(g_u))^0)=\uL.\endalign$$
Since $\boc=\cup_{l'\in\uL}l'g_s\fc l'{}\i$ we see that $y\boc y\i=\boc$. We have 
$ygy\i=yg_sy\i yg_uy\i=g_slg_ul\i=lgl\i\in\la g\ra$. Hence $y\la g\ra y\i=\la g\ra$. We
see that $y\in\cg'$. Thus, $\cg\sub\cg'$. The inclusion $\uL^F\sub\cg'$ is obvious. 
Hence $\uL^F\cg\sub\cg'$. Conversely, let $y\in\cg'$. Then $y=l'y'$ where $l'\in\uL^F$,
$y'\in\cg'$, $y'gy'{}\i=g$. We have $y'\uL y'{}\i=\uL$ hence $y'Ly'{}\i=L$. We have 
$y'\boc y'{}\i=\boc$, $y'(g_uL)y'{}\i=g_uL$ hence $y'\fc y'{}\i=\fc$. Thus $y'\in\cg$ 
and (a) holds.

\subhead 21.6\endsubhead
Let $(L,S)\in\AA$ and let $\ce\in\cs(S)$ be an irreducible cuspidal local system on $S$
(relative to $N_GL$). Let $Y=Y_{L,S},\tY=\tY_{L,S},\p:\tY@>>>Y$ be as in 3.13. Let 
$\tce$ be the local system on $\tY$ defined in 5.6. Let 
$\ti\G=\{n\in N_{G^0}L;nSn\i=S,Ad(n)^*\ce\cong\ce\}$ and let $\G=\ti\G/L$. Let $\bY$ be
the closure of $Y$ in $G$ and let $\fK=IC(\bY,\p_!\tce)$. Assume that $FL=L,FS=S$,
$F^*\ce\cong\ce$. Now $F:G@>>>G$ induces isomorphisms $F:\ti\G@>>>\ti\G$, $F:\G@>>>\G$.

For any $n\in\ti\G^F$ there is a well defined element $\et(n)\in\bbq^*$ such that the 
following holds: if $\a:\Ad(n)^*\ce@>>>\ce,\e:F^*\ce@>>>\ce$ are isomorphisms, then for
any $g\in S$, the composition $\ce_{nF(g)n\i}@>\a>>\ce_{F(g)}@>\e>>\ce_g$ is $\et(n)$ 
times the composition $\ce_{nF(g)n\i}@>\e>>\ce_{ngn\i}@>\a>>\ce_g$. (This follows from 
the irreducibility of $\ce$ and Schur's lemma.) Clearly, $\et(n)$ is independent of the
choice of $\a,\e$. It follows that, if $g\in S^F$ then the composition 
$\ce_{ngn\i}@>\a>>\ce_g@>\e>>\ce_g$ is $\et(n)$ times the composition 
$\ce_{ngn\i}@>\e>>\ce_{ngn\i}@>\a>>\ce_g$. Hence 
$\c_{\ce,\e}(g)=\et(n)\c_{\ce,\e}(ngn\i)$. This property characterizes $\et(n)$ since 
$S^F\ne\em$. Since $\c_{\ce,\e}$ is constant on $L^F$-conjugacy classes in $S^F$ it
follows that $\et$ induces a homomorphism $\ti\G^F/L^F=(\ti\G/L)^F=\G^F@>>>\bbq^*$. We
say that $(L,S,\ce)$ is {\it effective} if the associated homomorphism 
$\et:\G^F@>>>\bbq^*$ is identically $1$. (In this case we have 
$\c_{\ce,\e}(ngn\i)=\c_{\ce,\e}(g)$ for any $g\in S^F,n\in\ti\G^F$.)

We fix an isomorphism $\e:F^*\ce@>\si>>\ce$. Now $F$ induces 
Frobenius maps on $Y,\tY,\bY$. Also, $\e$ induces an isomorphism $F^*\tce@>\si>>\tce$ 
of local systems on $\tY$, an isomorphism $F^*\p_!\tce@>\si>>\p_!\tce$ of local systems
on $Y$ and an isomorphism $\ph:F^*\fK@>\si>>\fK$ in $\cd(\bY)$. As in 7.10, let $\EE$ 
be the algebra of endomorphisms of the local system $\p_!\tce$ on $Y$. As in 7.10(a) we
have a canonical decomposition $\EE=\op_{w\in\G}\EE_w$ where $\EE_w$ is a one 
dimensional subspace of $\EE$. From the definitions we see that for $w,y\in\G$ we have 
$\EE_w\EE_y=\EE_{wy}$. Let $b_w$ be a basis element of $\EE_w$. Then $\EE,\EE_w,b_w$ 
are as in 20.1. As in 20.3, let $V_i,(i\in[1,r'])$ be a set of representatives for the 
isomorphism classes of simple $\EE$-modules. We have canonically $\EE=\End(\fK)$. For 
$i\in[1,r']$ let $(\p_!\tce)_i=\Hom_\EE(V_i,\p_!\tce)$, $\fK_i=\Hom_\EE(V_i,\fK)$. Then
$(\p_!\tce)_i$ is an irreducible local system on $Y$ and $\fK_i=IC(\bY,(\p_!\tce)_i)$. 
Also, $(\p_!\tce)_i\not\cong(\p_!\tce)_{i'}$ for $i\ne i'$. We have canonically 
$\p_!\tce=\op_{i\in[1,r']}V_i\ot(\p_!\tce)_i$ and $\fK=\op_{i\in[1,r']}V_i\ot\fK_i$. 
For $f\in\Hom(\p_!\tce,\p_!\tce)$ we have
$$F^*(f)\in\Hom(F^*\p_!\tce,F^*\p_!\tce)=\Hom(\p_!F^*\tce,\p_!F^*\tce)=
\Hom(\p_!\tce,\p_!\tce)$$
where the last equality is obtained by using twice the isomorphism $F^*\tce@>>>\tce$ as
above. Hence we have a map $\io:\EE@>>>\EE,f\m F^*(f)$ which is an algebra isomorphism;
it carries $\EE_w$ onto $\EE_{F(w)}$ for any $w\in\G$. As in 20.3 we may assume that 
for $i\in[1,r]$, property 20.3($*$) holds and that for $i>r$ that property does not 
hold. For $i\in[1,r]$ we choose an isomorphism $\io_i:V_i@>>>V_i$ as in 20.3($*$). For 
any $w\in\G$, the isomorphism $b_w\ph:F^*\fK@>\si>>\fK$ corresponds under
$\fK=\op_{i\in[1,r']}V_i\ot\fK_i$ to an isomorphism

(a) $\op_{i\in[1,r']}V_i\ot F^*\fK_i@>>>\op_{i\in[1,r']}V_i\ot\fK_i$
\nl
which is an isomorphism of the summand $V_i\ot F^*\fK_i$ onto a summand 
$V_{i'}\ot\fK_{i'}$ where $i'=i$ for $i\in[1,r]$ and $i'\ne i$ for $i>r$; moreover, the
restriction of (a) to $V_i\ot F^*\fK_i,i\in[1,r]$, is of the form $b_w\io_i\ot\ph_i$ 
where $\ph_i:F^*\fK_i@>\si>>\fK_i$ is an isomorphism independent of $w$. (Note that 
$F^*\fK_i\not\cong\fK_i$ for $i>r$.) Taking induced maps on stalks and taking traces we
obtain for any $j\in\ZZ,g\in\bY^F$:
$$\tr(b_w\ph,\ch^j_g\fK)=\sum_{i\in[1,r]}\tr(b_w\io_i,V_i)\tr(\ph_i,\ch^j\fK_i).$$
Taking alternating sum over $j$ we obtain
$$\c_{\fK,b_w\ph}=\sum_{i\in[1,r]}\tr(b_w\io_i,V_i)\c_{\fK_i,\ph_i}.\tag b$$
We choose an element $g_w\in G^0$ such that $g_w\i F(g_w)=n_w$, a representative of $w$
in $\ti\G$. We set $L^w=g_wLg_w\i,S^w=g_wSg_w\i$, $\ce^w=\Ad(g_w\i)^*\ce$ (a local 
system on $S^w$). Then $F(L^w)=L^w$ and $F(S^w)=S^w$. We define an isomorphism
$\e^w:F^*\ce^w@>\si>>\ce^w$ in terms of $\e:F^*\ce@>\si>>\ce$ and $b_w$ as follows. By 
definition, $b_w$ defines for each $g\in S$ an isomorphism of stalks 
$\ce_{n_wgn_w\i}@>\si>>\ce_g$; hence it defines for any $g'\in S^w$ an isomorphism
$\ce_{n_wF(g_w)\i F(g')F(g_w)n_w\i}@>\si>>\ce_{F(g_w)\i F(g')F(g_w)}$ or equivalently 
$\ce_{g_w\i F(g')g_w}@>\si>>\ce_{F(g_w\i g'g_w)}$. Composing with 
$\e:\ce_{F(g_w\i g'g_w)}@>\si>>\ce_{g_w\i g'g_w}$ we obtain 
$\ce_{g_w\i F(g')g_w}@>\si>>\ce_{g_w\i g'g_w}$ that is, 
$\ce^w_{F(g')}@>\si>>\ce^w_{g'}$ which
comes from an isomorphism $\e^w:F^*\ce^w@>\si>>\ce^w$. We define $\p^w:\tY^w@>>>Y^w$,
$\tce^w,\fK^w$, $\ph^w:F^*\fK^w@>\si>>\fK^w$ in terms of $L^w,S^w,\ce^w,\e^w$ in the 
same way as $\p:\tY@>>>Y$, $\tce,\fK$, $\ph:F^*\fK@>\si>>\fK$ were defined in terms of 
$L,S,\ce,\e$. We have $Y^w=Y$ and the map $(g,xL)\m(g,xg_w\i L^w)$ is an isomorphism
$\mu:\tY@>>>\tY^w$ commuting with the projections $\p,\p^w$ onto $Y$. We have
$\mu^*\tce^w=\tce$ canonically. Hence $\mu$ induces an isomorphism
$\p_!\tce@>\si>>\p^w_!\tce^w$ hence an isomorphism $\mu':\fK@>\si>>\fK^w$. From the
definitions we see that the compositions $F^*\fK@>F^*\mu'>>F^*\fK^w@>\ph^w>>\fK^w$, 
$F^*\fK@>b_w\ph>>\fK@>\mu'>>\fK^w$ coincide. Hence for $j\in\ZZ,g\in\bY^F$ we have 
$\tr(b_w\ph,\ch^j_g\fK)=\tr(\ph^w,\ch^j_g\fK^w)$. Taking alternating sum over $j$ gives
$\c_{\fK,b_w\ph}=\c_{\fK^w,\ph^w}$. Introducing this in (b) we obtain
$$\c_{\fK^w,\ph^w}=\sum_{i\in[1,r]}\tr(b_w\io_i,V_i)\c_{\fK_i,\ph_i}.\tag c$$
Using (c) for $w$ running through $\baG$ (a set of representatives for the effective
$F$-twisted conjugacy classes in $\G$) and 20.4(e),(g), we see that

(d) {\it the functions $(\c_{\fK^w,\ph^w})_{w\in\baG}$ span the same vector space as 
the functions $(\c_{\fK_i,\ph_i})_{i\in[1,r]}$; moreover $|\baG|=r$.}
\nl
From the definitions we see that 

(e) {\it $(L^w,S^w,\ce^w)$ is effective if and only if $w\in\G$ is effective.}
\nl
We show that, 

(f) {\it if $x\in\G$ is not effective, then $\c_{\fK^x,\ph^x}=0$.}
\nl
It is enough to show that, for any $j\in\ZZ,g\in\bY^F$ we have 
$\tr(b_x\ph,\ch^j_g\fK)=0$. The proof is a repetition of that of 20.4(a). We can find 
$y\in\G_x$ such that $\g_x(y)\ne 1$. (Notation of 20.4.) We have
$$\align&\tr(\ph b_x,\ch^j_g\fK)=\tr(b_y\i\ph b_xb_y,\ch^j_g\fK)
=\g_x(y)\i\tr(b_y\i\ph\io\i(b_y)b_x,\ch^j_g\fK)\\&=\g_x(y)\i\tr(\ph b_x,\ch^j_g\fK).
\endalign$$
Thus $(1-\g_x(y)\i)\tr(\ph b_x,\ch^j_g\fK)=0$ and $\tr(\ph b_x,\ch^j_g\fK)=0$ as 
claimed.

\subhead 21.7\endsubhead
We preserve the setup of 21.6. Let $\Xi$ be the set of all triples $(L',S',[\ce'])$ 
such that

$(L',S')\in\AA$, 

$[\ce']$ is the isomorphism class of an irreducible cuspidal local system 
$\ce'\in\cs(S')$ (relative to $N_GL'$),

there exists $g\in G^0$ such that $gLg\i=L',gSg\i=S'$ and $\Ad(g\i)^*\ce\in[\ce']$.
\nl
Note that $G^0$ acts naturally on $\Xi$ and this action is transitive. The isotropy 
group of $(L,S,[\ce])$ in $G^0$ is $\ti\G$. Thus we may identify $\Xi=G^0/\ti\G$. Since
$\ti\G$ is $F$-stable, we have an induced Frobenius map $F:\Xi@>>>\Xi$ whose fixed 
point set consists of all $(L',S',[\ce'])\in\Xi$ such that $F(L')=L',F(S')=S'$ and
$F^*\ce'\cong\ce'$. By 19.7(a), the triples $(L^w,S^w,\ce^w)$ (where $w$ runs through a
set of representatives of the $F$-twisted conjugacy classes in $\G$) form a set of 
representatives for the $G^{0F}$-orbits in $\Xi^F$. Using now 21.6(e) we see that the 
triples $(L^w,S^w,\ce^w)$ (where $w$ runs through $\baG$) form a set of representatives
for the $G^{0F}$-orbits on the set of effective triples in $\Xi^F$. 

Let $({}^hL,{}^hS,[{}^h\ce])_{h\in I}$ be a set of representatives for the 
$G^{0F}$-orbits on the set of effective triples in $\Xi^F$. For each $h\in I$ choose 
${}^h\ce\in[{}^h\ce]$ and an isomorphism ${}^h\e:F^*({}^h\ce)@>\si>>{}^h\ce$. Define 
${}^h\fK$, ${}^h\ph:F^*({}^h\fK)@>\si>>{}^h\fK$ in terms of 
${}^hL,{}^hS,{}^h\ce,{}^h\e$ 
in the same way as $\fK$, $\ph:F^*\fK@>\si>>\fK$ were defined in terms of $L,S,\ce,\e$.
Let $(A_j)_{j\in J}$ be a set of representatives for the isomorphism classes of simple 
intersection cohomology complexes $A$ on $\bY$ that are summands of $\fK$ and satisfy 
$F^*A\cong A$. For each $j\in J$ choose an isomorphism $\ps_j:F^*A_j@>\si>>A_j$. We can
now reformulate 21.6(d) as follows.

(a) {\it the functions $(\c_{{}^h\fK,{}^h\ph})_{h\in I}$ span the same vector space as 
the functions \lb
$(\c_{A_j,\ps_j})_{j\in J}$; moreover, $|I|=|J|$.}

\subhead 21.8\endsubhead
Let $L,S,\ce$ be as in 21.6. Assume that $S$ contains a unipotent $L$-conjugacy class 
$\fc$ (necessarily unique hence $F$-stable), that $\ce$ is the inverse image under 
$S@>>>\fc,g\m g_u$ of an irreducible cuspidal local system $\cf$ on $\fc$ and that 
$\e:F^*\ce@>\si>>\ce$ is induced via $S@>>>\fc$ by an isomorphism 
$\e_0:F^*\cf@>\si>>\cf$. We show that

(a) {\it $(L,S,\ce)$ is effective.} 
\nl
Let $\bY,\fK,\ph:F^*\fK@>\si>>\fK,\EE,b_w,\G,F:\G@>>>\G,\io:\EE@>>>\EE$ be as in 21.6. 
By 21.6(e) it is enough to show that $1\in\G$ is effective (relative to the basis $b_w$
of $\EE$ and $F:\G@>>>\G$). We may assume that $b_w$ are chosen as in the proof of
Proposition 11.9. We must show that $y\in\G,F(y)=y\implies \io\i(b_y)=b_y$. Hence it is
enough to show that $\io(b_y)=b_{F(y)}$ for any $y\in\G$. Let $\boc$ be the unique 
unipotent $G^0$-conjugacy class of $G$ that is open dense in $\bY\cap G_{un}$ (see 
10.3). Let $\ch=\ch^0\fK|_\boc$, an irreducible local system on $G^0$ (see 11.8). The 
natural action of $\EE$ on $\fK$ induces an action of $\EE$ on $\ch$ in which $b_y$ 
acts as the identity map (by the choice of $b_y$). Let $f\in\Hom(\fK,\fK)=\EE$. The
commutative diagram
$$\CD
F^*\fK@>F^*f>>F^*\fK\\
@V\ph VV     @V\ph VV \\        
\fK  @>\io(f)>>\fK
\endCD$$
(which comes from the definition of $\io$) induces a commutative diagram
$$\CD
F^*\ch@>F^*f>>F^*\ch\\
@V\ph VV     @V\ph VV \\        
\ch  @>\io(f)>>\ch
\endCD$$
If $f=b_y$ then $f:\ch@>>>\ch$ is the identity map hence $F^*f:F^*\ch@>>>F^*\ch$ is the
identity map. Then the last commutative diagram shows that $\io(f):\ch@>>>\ch$ is the 
identity map. Since $\io(b_y)$ is a scalar multiple of $b_{F(y)}$ and $b_{F(y)}$ acts
on $\ch$ as the identity map, it follows that $\io(b_y)=b_{F(y)}$, as required.

\subhead 21.9\endsubhead
Let $\cv$ be the vector space of functions $G_{un}^F@>>>\bbq$ that are constant on 
$G^{0F}$-conjugacy classes in $G_{un}^F$. Let $\cn$ be the set of all pairs 
$(\boc,[\cf])$ where $\boc$ is a unipotent $G^0$-conjugacy class in $G$ and $[\cf]$ is 
an isomorphism class of an irreducible $G^0$-equivariant local system $\cf$ on $\boc$.
Define $F:\cn@>>>\cn$ by $F(\boc,[\cf])=(F(\boc),[F_!\cf])$. The fixed point set 
$\cn^F$ is the set of all $(\boc,[\cf])\in\cn$ such that $F(\boc)=\boc$ and
$F^*\cf\cong\cf$. For any $(\boc,[\cf])\in\cn^F$ we choose a local system 
$\cf\in[\cf]$ and an isomorphism $\ph_{\cf}:F^*\cf@>\si>>\cf$. The function 
$\c_{\cf,\ph_{\cf}}:\boc^F@>>>\bbq$ will be regarded as a function 
$G_{un}^{0F}@>>>\bbq$, equal to zero on $G_{un}^{0F}-\boc^F$. Using Lemma 19.7, we see 
that

(a) {\it for any $F$-stable unipotent $G^0$-conjugacy class $\boc'$ in $G$, the 
functions $\c_{\cf',\ph_{\cf'}}$ with $(\boc',[\cf'])\in\cn^F$ form a basis for the 
vector space of functions in $\cv$ that are zero on $G_{un}^F-\boc'{}^F$.}
\nl
From (a) we deduce

(b) {\it the functions $\c_{\cf,\ph_{\cf}}$ with $(\boc,[\cf])\in\cn^F$ form a basis 
of the vector space} $\cv$.
\nl
For $(\boc,[\cf])\in\cn^F$ let $\cf^\sh=IC(\bboc,\cf)$. Now $\ph_{\cf}$ induces an
isomorphism $\ph_{\cf}^\sh:F^*\cf^\sh@>\si>>\cf^\sh$ in $\cd(\bboc)$. Hence 
$\c_{\cf^\sh,\ph_{\cf}^\sh}:\bboc^F@>>>\bbq$ is well defined. We regard 
$\c_{\cf^\sh,\ph_{\cf}^\sh}$ as a function $G_{un}^F@>>>\bbq$, equal to zero on 
$G_{un}^F-\bboc^F$. This function is constant on $G^{0F}$-conjugacy classes. Hence it 
is of the form $\sum_{\boc'}c_{\boc'}f_{\boc'}$ where $\boc'$ runs over the unipotent 
$G^0$-conjugacy classes in $G$ such that $F(\boc')=\boc',\boc'\sub\bboc$,
$c_{\boc'}\in\bbq$ and $f_{\boc'}\in\cv$ is zero on $G^F-\boc'{}^F$. For such $\boc'$, 
$f_{\boc'}$ is a linear combination of functions $\c_{\cf',\ph_{\cf'}}$ with $\cf'$ 
such that $(\boc',[\cf'])\in\cn^F$ (see (a)). From the definitions we have 
$f_\boc=\c_{\cf,\ph_{\cf}}$. We see that
$$\c_{\cf^\sh,\ph_{\cf}^\sh}=\sum_{(\boc',[\cf'])\in\cn^F}
c_{[\cf],[\cf']}\c_{\cf',\ph_{\cf'}}$$
where $c_{[\cf],[\cf']}\in\bbq$ are uniquely determined and equal to zero unless 
$\boc'\sub\bboc$; moreover, if $\boc'=\boc$ then $c_{[\cf],[\cf']}=\d_{[\cf],[\cf']}$. 
Thus, the functions 
$\c_{\cf^\sh,\ph_{\cf}^\sh}$ are related to the functions $\c_{\cf,\ph_{\cf}}$ by an 
upper triangular matrix with all diagonal entries equal to $1$. Hence (b) implies:

(c) {\it the functions $\c_{\cf^\sh,\ph_{\cf}^\sh}$ with $(\boc,[\cf])\in\cn^F$ form 
a basis of the vector space} $\cv$.

\subhead 21.10\endsubhead
Let $\cy$ be the set of triples $(L,\fc,[\ff])$ where $L$ is the Levi of some parabolic
of $G^0$, $\fc$ is a unipotent $L$-conjugacy class of $N_GL$ with $\fc\sub N_G^\bul L$
and $[\ff]$ is the isomorphism class of an irreducible cuspidal local sytem $\ff$ on 
$\fc$ (relative to $N_GL$). Let $G^0\bsl\cy$ be the set of orbits of the 
natural $G^0$-action on $\cy$ given by conjugation all factors. Define $F:\cy@>>>\cy$ 
by $F(L,\fc,[\ff])=(F(L),F(\fc),[F_!\ff])$; we have $F(gy)=F(g)F(y)$ for all 
$g\in G^0$, $y\in\cy$. Hence $F$ induces a bijection $F:G^0\bsl\cy@>>>G^0\bsl\cy$. 
Putting together the generalized Springer correspondences 11.10(a) for the various 
connected components of $G$ that contain unipotent elements we obtain a canonical 
surjective map $\cn@>>>G^0\bsl\cy$. From the definitions we see that this map is 
compatible with the $F$-actions, hence it restricts to a surjective map 
$\cn^F@>>>(G^0\bsl\cy)^F$ whose fibres form a partition of $\cn^F$ into subsets 
$\cn^F_{\co}$ indexed by the $F$-stable $G^0$-orbits $\co$ on $\cy$. Using 21.9(c), we 
see that 

(a) $\cv=\op_\co\cv_\co$ 
\nl
($\co$ as above) where $\cv_\co$ is the subspace of $\cv$ with basis formed by the 
functions $\c_{\cf^\sh,\ph_{\cf}^\sh}$ with $(\boc,[\cf])\in\cn^F_\co$. Since $\co$ is 
a homogeneous space for the connected group $G^0$, it contains some $F$-fixed point 
$(L,\fc,[\ff])$. We have $F(L)=L,F(\fc)=\fc$. We choose $\ff\in[\ff]$. We have 
$F^*\ff\cong\ff$; we choose an isomorphism $\e_1:F^*\ff@>\si>>\ff$. Let $S$ be the 
stratum of $N_GL$ that contains $\fc$, let $\ce$ be the inverse of $\cf$ under 
$S@>>>\fc,g\m g_u$ and let $\e:F^*\ce@>\si>>\ce$ be the isomorphism induced by $\e_1$. 
Then $(L,S,\ce,\e)$ are as in 21.6. 

We will apply 21.7(a) in our case. Restricting the functions in 21.7(a) to 
$\bY^{\o F}=\bY\cap G_{un}^F$ ($\bY$ as in 21.6) we see that

(b) {\it the functions $\c_{{}^h\fK,{}^h\ph}|_{\bY^{\o F}}$ ($h\in I$) span the same 
vector space as the functions $\c_{A_j,\ps_j}|_{\bY^{\o F}}$ ($j\in J$); moreover, 
$|I|=|J|$.}
\nl
From the definition of generalized Springer correspondence we see 
that the functions $\c_{A_j,\ps_j}|_{\bY^{\o F}}$, ($j\in J$), extended by $0$ on 
$G_{un}^F-\bY^{\o F}$, are up to non-zero scalars the same as the functions
$\c_{\cf^\sh,\ph_{\cf}^\sh}$ with $(\boc,[\cf])\in\cn^F_{\co}$. In particular, they 
form a basis of the vector space $\cv_{\co}$. Using now (b), we see that the functions 
$\c_{{}^h\fK,{}^h\ph}|_{\bY^{\o F}}$, ($h\in I$), extended by $0$ on 
$G_{un}^F-\bY^{\o F}$ (or equivalently, the generalized Green functions
$Q_{G,{}^hL,{}^h\fc,{}^h\ff,{}^h\e_1}:G_{un}^F@>>>\bbq$, see below) form a basis of 
the vector space $\cv_{\co}$. Here ${}^h\fc$ is the set of unipotent elements in 
${}^hS$, ${}^h\ff={}^h\ce|_{{}^h\fc}$ and ${}^h\e_1$ is the restriction of ${}^h\e$ to 
${}^h\fc$. In our case 

(c) {\it all triples in $\Xi^F$ (see 21.7) are effective,} 
\nl
since all elements of $\G$ are effective (see 21.6(e) and 21.8(a)). Letting now $\co$ 
vary and using (a), we obtain the following result.

\proclaim{Proposition 21.11} The generalized Green functions $Q_{G,L,\fc,\ff,\e_1}$ 
(where $(L,\fc,[\ff])$ runs through a set of representatives for the $G^{0F}$-orbits on
$\cy^F$ and for each $(L,\fc,[\ff])$ in this set we choose $\ff\in[\ff]$ and 
$\e_1:F^*\ff@>\si>>\ff$) form a basis of the $\bbq$-vector space $\cv$ of functions 
$G_{un}^F@>>>\bbq$ that are constant on $G^{0F}$-conjugacy classes in $G_{un}^F$. 
\endproclaim

\subhead 21.12\endsubhead
Define $F:\fA_G@>>>\fA_G$ by $F(L,\fc)=(F(L),F(\fc))$. Let $(L,\fc)\in\fA_G^F$. Then
$\cc_L(\fc)$ is well defined (see 19.9).

(a) For any $n\in G^{0F}$ such that $nLn\i=L,n\fc n\i=\fc$ and any 
$f\in\cc_L(\fc),g\in\fc^F$ we have $f(ngn\i)=f(g)$.
\nl
Indeed, we may assume that $f=\c_{\ff,\e_1}$ where $\ff$ is an irreducible cuspidal
local system on $\fc$ and $\e_1:F^*\ff@>\si>>\ff$ is an isomorphism. In that case the 
result follows from 21.10(c), since we have automatically $\Ad(n)^*\ff\cong\ff$ (see 
11.7(a)).

We define a $\bbq$-linear map $\cc_L(\fc)@>>>\cv,f\m Q_{G,L,\fc}^f$ by the requirement 
that $Q_{G,L,\fc}^f=Q_{G,L,\fc,\ff,\e_1}$ for any $f=\c_{\ff,\e_1}$ as above. It is 
clear that this linear map is well defined. Let $J_{L,\fc}$ be the its image. Note that
$J_{L,\fc}$ depends only on the $G^{0F}$-orbit of $(L,\fc)$. We can reformulate 
Proposition 21.11 as follows.

(b) {\it For any $(L,\fc)\in\fA_G^F$, the linear map $\cc_L(\fc)@>>>J_{L,\fc}$ is an 
isomorphism. We have a direct sum decomposition $\cv=\op_{(L,\fc)}J_{L,\fc}$ where 
$(L,\fc)$ runs through a set of representatives for the $G^{0F}$-orbits on $\fA_G^F$.}
\nl
On $\op_{(L,\fc)\in\fA_G^F}\cc_L(\fc)$ we have a linear $G^F$-action: an element 
$g\in G^F$ takes $f\in\cc_L(\fc)$ to $f'\in\cc_{gLg\i}(g\fc g\i)$ where 
$f'(h)=f(g\i hg)$ for $h\in(g\fc g\i)^F$. This action restricts to a $G^{0F}$-action. 
Consider the linear map $\op_{(L,\fc)\in\fA_G^F}\cc_L(\fc)@>>>\cv$ whose restriction to
any summand $\cc_L(\fc)$ is $f\m Q_{G,L,\fc}^f$. Restricting this to the space of 
$G^{0F}$-invariants we obtain 

(c) {\it an isomorphism} $(\op_{(L,\fc)\in\fA_G^F}\cc_L(\fc))^{G^{0F}}@>\si>>\cv$. 
\nl
This follows immediately from (b),(a).

\subhead 21.13\endsubhead
Let $s\in G^F$ be a semisimple element and let $(L,\fc)\in\fA_{Z_G(s)}^F$. Let
$\cc^s_L(\fc)$ be the vector space $\cc_L(\fc)$ defined as in 19.9 with respect to 
$Z_G(s)$ instead of $G$. Let ${}'\cc^s_L(\fc)$ be the subspace of $\cc^s_L(\fc)$ 
consisting of all functions that are invariant under the natural action of 
$\{g\in Z_{G^0}(s)^F;gLg\i=L,g\fc g\i=\fc\}$. Replacing $G$ by $Z_G(s)$ in 21.12(c) we 
obtain an isomorphism

$(\op_{(L,\fc)\in\fA_{Z_G(s)}^F}\cc^s_L(\fc))^{Z_G(s)^{0F}}@>\si>>\cv_s$
\nl
where $\cv_s$ is the vector space of functions 
$\{\text{unipotent elements in }Z_G(s)^F\}@>>>\bbq$ that are constant on 
$Z_G(s)^{0F}$-conjugacy classes. Taking now invariants for the natural action of 
$Z_{G^0}(s)^F$ (which contains $Z_G(s)^{0F}$ as a normal subgroup) we obtain an 
isomorphism

$(\op_{(L,\fc)\in\fA_{Z_G(s)}^F}\cc^s_L(\fc))^{Z_{G^0}(s)^F}@>\si>>\cv'_s$
\nl
or, equivalently, an isomorphism

($*$) $(\op_{(L,\fc)\in\fA_{Z_G(s)}^F}{}'\cc^s_L(\fc))^{Z_{G^0}(s)^F}@>\si>>\cv'_s$
\nl
where $\cv'_s$ is the vector space of functions 
$\{\text{unipotent elements in }Z_G(s)^F\}@>>>\bbq$ that are constant on 
$Z_{G^0}(s)^F$-conjugacy classes. We now take the direct sum of these isomorphisms over
all semisimple $s$ in $G^F$ and then take invariants for the natural action of 
$G^{0F}$. We obtain an isomorphism
$$(\op_{(s,L,\fc)\in\cx}{}'\cc^s_L(\fc))^{G^{0F}}@>>>VV\tag a$$
where $\cx$ is the set of all triples $(s,L,\fc)$ with $s\in G^F$ semisimple and 
$(L,\fc)\in\fA_{Z_G(s)}^F$ and $\VV$ is the vector space of functions $G^F@>>>\bbq$ 
that are constant on $G^{0F}$-conjugacy classes. 

Let $\cx^1$ be the set of all quadruples $(s,u,L,\fc)$ where $s\in G^F$ is semisimple,
$u\in G^F$ is unipotent, quasi-semisimple in $N_{Z_G(s)}L$, $(L,\fc)\in\fA_{Z_G(s)}^F$,
$\fc\sub uL$. Now $G^{0F}$ acts naturally on $\cx^1$ and the map $\cx^1@>>>\cx$, 
$(s,u,L,\fc)\m(s,L,\fc)$ induces a bijection 

(b) $G^{0F}\bsl\cx^1@>\si>>G^{0F}\bsl\cx$
\nl
on the sets of $G^{0F}$-orbits. (We use the fact that, for fixed $(s,L,\fc)\in\cx$, the
set of unipotent quasi-semisimple elements of $N_{Z_G(s)}L$ that are fixed by $F$ and 
are contained in the component $\fc L$ of $N_{Z_G(s)}L$ is a single $L^F$-conjugacy 
class; this follows from 19.11.) Now $(g,L,\fc)\m(g_s,g_u,L,\fc)$ is a bijection

(c) $\car'@>\si>>\cx^1$.
\nl
(We use 1.4(c) and 21.2(a).) Combining (a),(b),(c) we obtain an isomorphism
$$(\op_{(g,L,\fc)\in\car'}{}'\cc^{g_s}_L(\fc))^{G^{0F}}@>\si>>\VV.\tag d$$
Assume that $(g,\uL,\boc)\in\car$ corresponds to $(g,L,\fc)\in\car'$ under 21.4(c) and 
let $\la g\ra$ be the $\uL^F$-conjugacy of $g$. From 19.15(b),(c) applied to 
$N_G\uL,\uL$ instead of $G,G^0$ we have an isomorphism
$$\cc_{\uL,\la g\ra}(\boc)@>\si>>\cc_{H_{\uL}(g)}(\fc)\tag e$$
where $H_{\uL}(g)=\{l\in\uL; lg_sl\i=g_s, lg_ul\i\in Z_{\uL}(g_s)^0g_u\}$.

Let $\cg=\{y\in Z_{G^0}(g_s)^F;yLy\i=L,y\fc y\i=\fc\}$ (a group containing 
$H_{\uL}(g)^F$ as a normal subgroup).

Let $\cg'=\{y\in G^{0F},y\uL y\i=\uL,y\boc y\i=\boc,y\la g\ra y\i=\la g\ra\}$ (a group
containing $\uL^F$ as a normal subgroup).

Assume that $f\m\baf$ under (e). We show that the following two conditions are 
equivalent:

(i) $\baf$ is invariant under the natural action of $\cg$;

(ii) $f$ is invariant under the natural action of $\cg'$.
\nl
Assume that (i) holds. Let $y\in\cg'$. Then $y=l'y'$ where $l'\in\uL^F,y'\in\cg$ (see
21.5(a)). We must show that $f(ylg_sul\i y\i)=f(lg_sul\i)$ for $l\in\uL^F,u\in\fc^F$ or
that $f(l'y'lg_sul\i y'{}\i l'{}\i)=\baf(u)$ or that $\baf(y'uy'{}\i)=\baf(u)$; this 
follows from $y'\in\cg$.

Assume that (ii) holds. Let $y\in\cg$. Then $y\in\cg'$ (see 21.5(a)). We must show that
$\baf(yuy\i)=\baf(u)$ for $u\in\fc^F$ or that $f(g_syuy\i)=f(g_su)$ or that 
$f(yg_suy\i)=f(g_su)$; this follows from $y\in\cg'$.

From the equivalence of (i),(ii), we see that (e) restricts to an isomorphism
$${}'\cc_{\uL,\la g\ra}(\boc)@>\si>>{}'\cc^{g_s}_L(\fc)\tag f$$
where ${}'\cc_{\uL,\la g\ra}(\boc)$ is the space of all functions 
$f\in\cc_{\uL,\la g\ra}(\boc)$ that are invariant under the natural action of 
$\{y\in G^{0F},y\uL y\i=\uL,y\boc y\i=\boc,y\la g\ra y\i=\la g\ra\}$. Using (f), we 
deduce from (d) an isomorphism
$$(\op_{(g,\uL,\boc)\in\car}{}'\cc_{\uL,\la g\ra}(\boc))^{G^{0F}}@>\si>>\VV$$
or equivalently an isomorphism
$$(\op_{(\uL,\boc)\in\fA_G^F}\op_\g{}'\cc_{\uL,\g}(\boc))^{G^{0F}}@>>>\VV$$
where $\g$ runs through the set of $\uL^F$-orbits on $(\s_{N_G\uL}\boc)^F$, or 
equivalently an isomorphism
$$(\op_{(\uL,\boc)\in\fA_G^F}\op_\g\cc_{\uL,\g}(\boc))^{G^{0F}}@>\si>>\VV.$$
From the definitions we have canonically $\op_\g\cc_{\uL,\g}(\boc)=\cc_{\uL}(\boc)$ 
hence we obtain an isomorphism
$$(\op_{(\uL,\boc)\in\fA_G^F}\cc_{\uL}(\boc))^{G^{0F}}@>\si>>\VV.\tag g$$
We define $F:\AA@>>>\AA$ by $F(\uL,S)=(F(\uL),F(S))$. There is a well defined 
surjective map $\fA_G^F@>>>\AA^F$ given by $(\uL,\boc)\m(\uL,S)$ where $S$ is defined
by $\boc\sub S$. Moreover, if $(\uL,S)\in\AA^F$ is given we have a natural isomorphism
$\cc_{\uL}(S)@>\si>>\op_\boc\cc_{\uL}(\boc)$ where $\boc$ runs over the $F$-stable
$\uL$-conjugacy classes contained in $S$. (A special case of 19.10(b).) Introducing
this in (g) we obtain an isomorphism
$$(\op_{(\uL,S)\in\AA^F}\cc_{\uL}(S))^{G^{0F}}@>\si>>\VV.$$
For each $(\uL,S)\in\AA^F$ we have a canonical direct sum decomposition
$$\cc_{\uL}(S)=\op_{[\ce]}\cc_{\uL}^{[\ce]}(S)$$
where $[\ce]$ runs over the set of isomorphism classes of irreducible cuspidal local 
systems $\ce\in\cs(S)$ (relative to $N_G\uL$) such that $F^*\ce\cong\ce$ and
$\cc_{\uL}^{[\ce]}(S)$ is the line spanned by $\c_{\ce,\e}$ where $\ce\in[\ce]$ and
$\e:F^*\ce@>\si>>\ce$. (This follows from 19.8(a).) Hence we have an isomorphism
$$(\op_{(\uL,S,[\ce])\in\ti\AA^F}\cc_{\uL}^{[\ce]}(S))^{G^{0F}}@>\si>>\VV\tag h$$
where $\ti\AA^F$ is the set of triples $(\uL,S,[\ce])$ with $(\uL,S)\in\AA^F$ and 
$[\ce]$ is as above. The left hand side of (h) is naturally a direct sum of subspaces 
$V_{(\uL,S,[\ce])}$ indexed by a set of representatives for the $G^{0F}$-orbits on 
$\ti\AA^F$ and $V_{(\uL,S,[\ce])}$ is the space of vectors in the one dimensional 
vector space $\cc_{\uL}^{[\ce]}(S)$ that are invariant under the natural action of the 
group $\{g\in G^{0F};gLg\i=L,gSg\i=S,\Ad(g)^*\ce\cong\ce\}$. From the definitions we 
see that $V_{(\uL,S,[\ce])}$ is $1$-dimensional if $(\uL,S,\ce)$ is effective and is 
$0$ if $(\uL,S,\ce)$ is not effective. 

Thus the left hand side of (h) has dimension equal to the number of $G^{0F}$-orbits on 
the set of effective triples in $\ti\AA^F$. Hence this number is equal to the dimension
of the right hand side that is, to the number of $G^{0F}$-conjugacy classes in $G^F$.

\proclaim{Theorem 21.14}Let $\ca$ be a set of representatives for the $G^{0F}$-orbits 
on the set of effective triples in $\ti\AA^F$. For each $(\uL,S,[\ce])\in\ca$ we choose
$\ce\in[\ce]$ and an isomorphism $\e:F^*\ce@>\si>>\ce$. To $\uL,S,\ce,\e$ we associate 
$\fK\in\cd(G)$ and $\ph:F^*\fK@>\si>>\fK$ as in 21.6 (with $\uL$ instead of $L$). The 
functions $\c_{\fK,\ph}$ (one for each $(\uL,S,[\ce])\in\ca$) form a $\bbq$-basis of
the vector space $\VV$ of functions $G^F@>>>\bbq$ that are constant on 
$G^{0F}$-conjugacy classes.
\endproclaim
Let $\VV'$ be the subspace of $\VV$ spanned by the functions $\c_{\fK,\ph}$ in the 
theorem. By the last paragraph in 21.13, it is enough to show that $\VV'=\VV$. This 
will be done in 21.17. Note that in the definition of $\VV'$ we may include the
functions $\c_{\fK,\ph}$ corresponding to non-effective triples in $\ti\AA$ (these
functions are identically $0$ by 21.6(f)).

\subhead 21.15 \endsubhead
Let $(L,S)\in\AA^F$. We define a linear function $\Ps:\cc_L(S)@>>>\VV$ by the 
requirement that for any irreducible cuspidal local system $\ce\in\cs(S)$ and any 
$\e:F^*\ce@>\si>>\ce$ we have $\Ps(\c_{\ce,\e})=\c_{\fK,\ph}$ where $\fK,\ph$ are
defined as in 21.6. From Theorem 16.14 we see that for any $f\in\cc_L(S)$ and any 
$y\in G^F$ we have
$$\Ps(f)(y)=\sum_{x\in G^{0F},\dd;x\i y_sx\in S_s}|L_x^F||Z_G(y_s)^{0F}|\i|L^F|\i
                                   Q_{L_x,Z_G(y_s),\dd}^{f_x^\dd}(y_u)\tag a$$
where $\dd$ runs over the set of $F$-stable $Z_G(y_s)^0$-orbits contained in

$\{v\in Z_G(y_s);v\text{ unipotent, }x\i y_svx\in S\}$, 
\nl
$L_x=Z_{xLx\i}(y_s)^0$ (a Levi of some parabolic of $Z_G(y_s)^0$), and 
$f_x^\dd\in\cc_{L_x}(\dd)$ is defined by $f_x^\dd(v)=f(x\i y_svx)$. (Notation of 
21.12.)

\subhead 21.16\endsubhead
Assume that we are given an $F$-stable $L$-conjugacy class $\boc$ in $S$ and a 
quasi-semisimple element $g\in N_GL$ such that $F(g)=g$ and such that

$\fc=\{u\in Z_L(s)^0g_u; u\text{ unipotent },su\in\boc\}\ne\em$.
\nl
Here $s=g_s$. We assume that $\boc$ is cuspidal (relative to $N_GL$). Then $\fc$ is a 
single $Z_L(g_s)^0$-conjugacy class (see 17.13). Assume that $f$ in 21.15 satisfies:

(i) $f|_{S^F-\boc^F}=0$,

(ii) $h\in\boc^F,f(h)\ne 0\implies h=lsul\i$ for some $l\in L^F$ and some $u\in\fc^F$.
\nl
Consider the function $\tf\in\VV$ given on $y\in G^F$ by

(a) $\tf(y)=0$ if $y_s$ is not $G^{0F}$-conjugate to $s$;
    
$\tf(y)=|Z_G(s)^{0F}|\i\sum_{z\in Z_{G^0}(s)^F}Q_{L_1,Z_G(s),\fc}^h(z\i y_uz)$ 
if $y_s=s$;
\nl
here $L_1=Z_L(y_s)^0$ and $h:\fc^F@>>>\bbq$ is given by $h(v)=f(y_sv)$. We show that 

(b) $\Ps(f)=|L_1^F||Z_L(s)^F|\i\tf$.
\nl
Consider the sum over $x,\dd$ in 21.15(a) for our $f$. If $f_x^\dd$ is not identically
$0$ then $f(x\i y_svx)\ne 0$ for some $v\in\dd^F$. Hence there exist $l\in L^F$,
$u\in\fc^F$ such that $x\i y_svx=lsul\i$ where $v\in\dd^F$. Thus,
$xls(xl)\i=y_s,xlu(xl)\i=v$. We see that $\Ps(f)(y)=0$ if $y_s$ is not 
$G^{0F}$-conjugate to $s$. Assume now that $y_s=s$. Setting $z=xl$ we have 
$z\in Z_{G^0}(s)^F,zuz\i=v$. Hence $\dd=z\fc z\i$. We see that
$$\align&\Ps(f)(y)=\sum_{z\in Z_{G^0}(s)^F,l\in L^F}|L_z^F||Z_G(s)^{0F}|\i|L^F|\i
|Z_L(s)^F|\i Q_{L_z,Z_G(s),z\fc z\i}^{f_z^{z\fc z\i}}(y_u)\\&
=|L_1^F||Z_G(s)^{0F}|\i|Z_L(s)^F|\i
\sum_{z\in Z_{G^0}(s)^F}
Q_{L_1,Z_G(s),\fc}^{f_1^{\fc}}(z\i y_uz),\endalign$$
as desired.

\subhead 21.17\endsubhead
From the definitions we see that $\VV'$ is the subspace of $\VV$ spanned by the
functions $\Ps(f)$ for various $(L,S)\in\AA^F$ and $f\in\cc_L(S)$. Since the vector
space $\cc_L(S)$ is spanned by functions as in 21.16(i),(ii) (see 19.10(b),19.14(c))
we see that any function $\tf$ as in 21.16(a) is contained in $\VV'$. In formula 
21.16(a) defining $\tf$ we may take $(g,L_1,\fc)$ to be any triple in $\car'$ and $h$ 
to be any function 
in $\cc_{L_1}(\fc)$ that is invariant under the natural action of 
$\{g'\in Z_{G^0}(s)^F;g'L_1g'{}\i=L_1,g'\fc g'{}\i=\fc\}$. (This follows from the 
bijection $\car\lra\car'$ in 21.4(c), the isomorphism 19.15(c) and the equivalence of 
(i),(ii) in 21.13.) 
But such $\tf$ span the vector space $\VV_s$ of functions in $\VV$ that vanish at 
elements whose semisimple parts are not $G^{0F}$-conjugate to $s$. ($\VV_s$ may be 
identified with $\cv'_s$ in 21.13 and we may use the isomorphism 21.13($*$).) We see 
that $\VV_s\sub\VV'$. Since $\VV$ is the sum of its subspaces $\VV_s$ for various $s$, 
we see that $\VV\sub\VV'$. Hence $\VV=\VV'$. Theorem 21.14 is proved. 

\subhead 21.18\endsubhead
Let $J$ be the set of all triples $(L,S,[\ce])$ where $(L,S)\in\AA$ and $[\ce]$ is the
isomorphism class of an irreducible cuspidal local system $\ce\in\cs(S)$. The group
$G^0$ acts on $J$ by $g:(L,S,[\ce])\m(\Ad(g)L,\Ad(g)S,[\Ad(g\i)^*\ce])$; let 
$G^0\bsl J$ be the set of orbits. 

Let $\fA(G)$ be the subcategory of $\cd(G)$ whose objects are the complexes $X$ on $G$
such that $X[d]$ is a (simple) admissible perverse sheaf on $G$ with support of 
dimension $d$ (see 6.7). 

Let $\ufA(G)$ be the set of isomorphism classes of objects in $\fA(G)$. We define a map
$$j:\ufA(G)@>>>G^0\bsl J$$
as follows. Let $A\in\fA(G)$. By definition there exists $(L,S,\ce)$ as above such that
$A$ is isomorphic to a direct summand of $IC(\bY,\p_!\tce)$ (extended by $0$ on 
$G-\bY$), with $\p,\tce,\bY$ as in 5.6. Then $j$ takes the isomorphism class of $A$ to 
the $G^0$-orbit of $(L,S,[\ce])$. To show that this is well defined we must show that, 
if $(L',S',\ce')$ is another triple like $(L,S,\ce)$ such that $A$ is isomorphic to a 
direct summand of $IC(\bY',\p'_!\tce')$ (extended by $0$ on $G-\bY'$), with 
$\p',\tce',\bY'$ defined as in 5.6 in terms of $L',S',\ce'$ instead of $L,S,\ce$, then
$(L,S,[\ce])$, $(L',S',[\ce'])$ are in the same $G^0$-orbit. Now $A$ is an intersection
cohomology complex supported by $\bY$ and also by $\bY'$. It follows that $\bY=\bY'$ 
hence $Y_{L,S}=Y_{L',S'}$. Using 3.12(b) we deduce that $(L,S),(L',S')$ are in the same
$G^0$-orbit. Hence we may assume that $L=L',S=S'$. Then $\ce,\ce'\in\cs(S)$ and there 
exists an irreducible local system on $\bY$ which is a direct summand of both 
$\p_!\tce$ and $\p_!\tce'$. Thus, $\Hom(\p_!\tce,\p_!\tce')\ne 0$. We now repeat an 
argument in 7.10 (and use notation there):
$$\align&\Hom(\p_!\tce,\p_!\tce')=\Hom(\p^*\p_!\tce,\tce')
=\op_{w\in\cw_S}\Hom(f_w^*\tce,\tce')\\&=\op_{w\in\cw_S}\Hom(a^*f_w^*\tce,a^*\tce')=
\op_{w\in\cw_S}\Hom(\hf_w^*a^*\tce,a^*\tce')\\&
=\op_{w\in\cw_S}\Hom(\hf_w^*b^*\ce,b^*\ce')
=\op_{w\in\cw_S}\Hom(b^*\Ad(n_w)^*\ce,b^*\ce')\\&=\op_{w\in\cw_S}
\Hom(\Ad(n_w)^*\ce,\ce').\endalign$$
We see that $\op_{w\in\cw_S}\Hom(\Ad(n_w)^*\ce,\ce')\ne 0$. Hence 
$\Hom(\Ad(n_w)^*\ce,\ce')\ne 0$ for some $w\in\cw_S$, so that $\Ad(n_w)^*\ce\cong\ce'$.
Thus, $(L,S,[\ce])$ is in the same $G^0$-orbit as $(L,S,[\ce'])$, as required.

\subhead 21.19\endsubhead
From the definitions we see that the map $j$ in 21.18 is compatible with the maps 
$F:\ufA(G)@>>>\ufA(G)$, $F:G^0\bsl J@>>>G^0\bsl J$ defined by $A\m F_!A$, 
$(L,S,[\ce])\m(F(L),F(S),[F_!\ce])$. Hence it induces a map 
$j_0:\ufA(G)^F@>>>(G^0\bsl J)^F$ on the fixed point sets of $F$. 

\subhead 21.20\endsubhead
From Theorem 21.14 we see that $\VV=\op_\Xi\VV^\Xi$ where $\Xi$ runs over the set
of $F$-stable $G^0$-orbits in $J$ and $\VV^\Xi$ is the subspace of $\VV$ with basis
given by the characteristic functions of $F$-stable effective triples 
$(L,S,[\ce])\in\Xi$ (up to the action of $G^{0F}$). Using now 21.7(a) we see that
another basis of $\VV^\Xi$ is given by the characteristic functions of objects in
$j_0\i(\Xi)$. (Either of these bases is defined only up to multiplication of any of its
members by a non-zero scalar.) Here we have used the fact that $\Xi$ contains at least
one $F$-stable triple $(L,S,[\ce])$ which follows from Lang's theorem for $G^0$ since 
$\Xi$ is a homogeneous $G^0$-space. Since $\ufA(G)^F=\sqc_\Xi j_0\i(\Xi)$, we see that
the following result holds.

\proclaim{Theorem 21.21} Let $\ca'$ be a set of representatives for the isomorphism 
classes of objects $A\in\fA$ such that $F^*A\cong A$. For each $A\in\ca'$ we choose an 
isomorphism $\a:F^*A@>\si>>A$. The functions $\c_{A,\a}$ (one for each $A\in\ca'$) form
a $\bbq$-basis of the vector space $\VV$ of functions $G^F@>>>\bbq$ that are constant 
on $G^{0F}$-conjugacy classes.
\endproclaim

\head 22. Twisted induction of class functions \endhead
\subhead 22.1\endsubhead
This section gives an application of Theorem 21.14 to the construction of a "twisted 
induction" map (see 22.3) from certain functions on a subgroup of $G^F$ to 
functions on $G^F$. 

\proclaim{Lemma 22.2} Let $L$ be a Levi of a parabolic of $G^0$ and let $L'$ be a Levi 
of a parabolic of $L$. Let $\d'$ be a connected component of $N_{N_GL}L'$ and let $\d$ 
be the connected component of $N_GL$ that contains $\d'$. Assume that 
$\d'\sub N_{N_GL}^\bul(L')$ and $\d\sub N_G^\bul L$. Then $\d'\sub N_G^\bul(L')$.
\endproclaim
Since $\d\sub N_G^\bul L$, there exists a parabolic $P$ of $G^0$ such that $L$ is a 
Levi of $P$ and $\d\sub N_GP$. Since $\d'\sub N_{N_GL}^\bul(L')$, there exists a 
parabolic $Q$ of $L$ such that $L'$ is a Levi of $Q$ and $\d'\sub N_{N_GL}Q$. Then 
$P'=QU_P$ is a parabolic of $G^0$ such that $L'$ is a Levi of $P'$. If $g\in\d'$ then 
$gQg\i=Q$ and $gU_Pg\i=U_P$ (since $g\in\d\sub N_GP$) hence $gP'g\i=P'$. Thus 
$\d'\sub N_GP'$. We see that $\d'\sub N_G^\bul(L')$, as required.

\subhead 22.3\endsubhead
Let $L$ be a Levi of a parabolic of $G^0$ and let $\d$ be a connected component of 
$N_GL$ contained in $N_G^\bul L$. We assume that $F(L)=L,F(\d)=\d$. Let $D$ be the 
connected component of $G$ that contains $\d$. Let $\VV_L(\d)$ (resp. 
$\VV_{G^0}(D)$) be the set of all functions $\d^F@>>>\bbq$ (resp. $D^F@>>>\bbq$) that 
are constant on $L^F$-conjugacy classes in $\d$ (resp. on $G^{0F}$-conjugacy classes in
$D$). 

There is a unique $\bbq$-linear map
$$R_\d^D:\VV_L(\d)@>>>\VV_{G^0}(D)$$
such that the following holds.

Let $L'$ be a Levi of a parabolic of $L$ with $F(L')=L'$, let $S'$ be an isolated
stratum of $N_{N_GL}(L')=N_GL\cap N_GL'$ with $F(S')=S',S'\sub\d$,
$S'\sub N_{N_GL}^\bul(L')$, let $\ce'$ be an irreducible cuspidal local system in 
$\cs(S')$ and let $\e':F^*\ce'@>\si>>\ce'$ be an isomorphism. Define 
$\fK'\in\cd(N_GL)$, $\ph':F^*\fK'@>\si>>\fK'$ in terms of $N_GL,L',S',\ce',\e'$ and 
$\fK''\in\cd(G)$, $\ph'':F^*\fK''@>\si>>\fK''$ in terms of 
$G,L',S',\ce',\e'$ in the same way as $\fK\in\cd(G),\ph:F^*\fK@>\si>>\fK$ were defined 
in 21.6 in terms of $G,L,S,\ce,\e$. (Note that $\fK''$ is well defined since 
$S'\sub N_G^\bul L'$, by Lemma 22.2.) Then 
$$R_\d^D(\c_{\fK',\ph'}|_{\d^F})=\c_{\fK'',\ph''}|_{D^F}.$$
To see that this definition is correct we use the fact that the functions
$\c_{\fK',\ph'}|_{\d^F}$ as above (with $\e'$ chosen for each $L',S',\ce'$ given up to
$L^F$-conjugacy) provide a basis for $\VV_L(\d)$ (which follows from Theorem 21.14 for 
$N_GL$ instead of $G$); note that the choice of $\e'$ is immaterial since the same 
choice is made in the definition of $\c_{\fK'',\ph''}|_{D^F}$.

\subhead 22.4\endsubhead
For any $G^{0F}$-conjugacy class $c$ of semisimple elements in $G^F$ let 
$\VV_{G^0,c}(D)$ be the vector space consisting of all functions in $\VV_{G^0}(D)$ 
that vanish on elements $g\in D^F$ with $g_s\notin c$. We have a direct sum 
decomposition 

(a) $\VV_{G^0}(D)=\op_{c}\VV_{G^0,c}(D)$ 
\nl
where $c$ runs over the semisimple $G^{0F}$-conjugacy classes in $G^F$. Similarly we 
have a direct sum decomposition 

(b) $\VV_L(\d)=\op_{c'}\VV_{L,c'}(\d)$
\nl
where $c'$ runs over the semisimple $L^F$-conjugacy classes in $(N_GL)^F$. The next
result shows that $R_\d^D$ is compatible with the direct sum decompositions (a),(b).

\proclaim{Proposition 22.5} Let $c'$ be any semisimple $L^F$-conjugacy class in 
$(N_GL)^F$ and let $c$ be the semisimple $G^{0F}$-conjugacy class in $G^F$ such that 
$c'\sub c$. Then $R_\d^D(\VV_{L,c'}(\d))\sub\VV_{G^0,c}(D)$.
\endproclaim
Let $L',S'$ be as in 22.3. As in 21.15 we define linear maps 
$\Ps':\cc_{L'}(S')@>>>\VV_L(\d)$, $\Ps'':\cc_{L'}(S')@>>>\VV_{G^0}(D)$ by the 
requirement that for any $\ce',\e'$ as in 22.3 we have 
$\Ps'(\c_{\ce',\e'})=\c_{\fK',\ph'}|_{\d^F}$,
$\Ps''(\c_{\ce',\e'})=\c_{\fK'',\ph''}|_{D^F}$ (notation of 22.3). Clearly,

(a) $R_\d^D(\Ps'(f))=\Ps''(f)$ 
\nl
for $f=\c_{\ce',\e'}$ hence also for any $f\in\cc_{L'}(S')$. Assume now that $\boc$ is 
an $F$-stable $L'$-conjugacy class in $S'$ and that $g$ is a quasi-semisimple element 
in $N_{N_GL}(L')$ such that $g\in\d,F(g)=g$ and such that
$\fc=\{u\in Z_{L'}(s)^0g_u; u\text{ unipotent },su\in\boc\}\ne\em$ (with $s=g_s$). 
Assume that $\boc$ is cuspidal (relative to $N_{N_GL}(L')$) so that $\fc$ is a single 
$Z_{L'}(s)^0$-conjugacy class. Assume that $f\in\cc_{L'}(S')$ satisfies

(i) $f|_{S'{}^F-\boc^F}=0$,

(ii) $h\in\boc^F,f(h)\ne 0\implies h=lsul\i$ for some $l\in L'{}^F$ and some
$u\in\fc^F$.
\nl
Consider the functions $\tf\in\VV_L(\d),\tf''\in\VV_{G^0}(D)$ defined by

$\tf(y)=0$ if $y\in\d^F$, $y_s$ is not $L^F$-conjugate to $s$;
    
$\tf(y)=|Z_L(s)^{0F}|\i\sum_{z\in Z_L(s)^F}Q_{L'_1,Z_{N_GL}(s),\fc}^h(z\i y_uz)$ 
if $y\in\d^F,y_s=s$;

$\tf''(y)=0$ if $y\in D^F$, $y_s$ is not $G^{0F}$-conjugate to $s$;
    
$\tf''(y)=|Z_G(s)^{0F}|\i\sum_{z\in Z_{G^0}(s)^F}Q_{L'_1,Z_G(s),\fc}^h(z\i y_uz)$ 
if $y\in D^F,y_s=s$;

here $L'_1=Z_{L'}(s)^0$ and $h:\fc^F@>>>\bbq$ is given by $h(v)=f(y_sv)$. Using 
21.16(f) for $G$ and for $N_GL$ and (a) we see that

(b) $R_\d^D(\tf)=\tf''$.
\nl
As in 21.17, here we may assume that $(g,L'_1,\fc)$ used in the definition of 
$\tf,\tf''$ is any triple in $\car'$ (for $N_GL$ instead of $G$) with $g\in\d$ and $h$ 
is any function in $\cc_{L'_1}(\fc)$ that is invariant under the natural action of 
$\{g'\in Z_L(s)^F;g'L'_1g'{}\i=L'_1,g'\fc g'{}\i=\fc\}$. Since the functions $\tf$ as 
above span the vector space $\VV_{L,c'}(\d)$ where $c'$ is the $L^F$-conjugacy class of
$g_s$ and the corresponding functions $\tf''$ are contained in the corresponding
$\VV_{G^0,c}(D)$ we see that the Proposition follows from (b).

\subhead 22.6\endsubhead
Let $L'\sub L\sub G^0,\d'\sub\d$ be as in Lemma 22.2 and let $D$ be the connected 
component of $G$ that contains $\d$. Then $R_{\d'}^\d,R_\d^D, R_{\d'}^D$ are well 
defined (the last one, by Lemma 22.2). From the definitions we see that the 
transitivity property

(a) $R_\d^D\circ R_{\d'}^\d=R_{\d'}^D$
\nl
holds.

\enddocument